\documentclass[12pt,reqno,a4paper]{amsart}
\usepackage{amssymb,amsfonts,latexsym,mathrsfs, tikz}

\title[Affine reflection subgroups of Coxeter groups]{Affine reflection subgroups of Coxeter groups}

\author{Xiang Fu}
\author{Lawrence Reeves}
\author{Linxiao Xu}

\dedicatory{\upshape
Beijing International Center for Mathematical Research\\
Peking University, Beijing, China\\
\texttt{fuxiang@math.pku.edu.cn}\\[1em]
School of Mathematics and Statistics \\
University of Melbourne, VIC 3010, Australia\\
\texttt{lreeves@unimelb.edu.au}\\[1em]
Department of Foundational Mathematics\\
School of Mathematics and Physics\\
Xi'an Jiaotong-Liverpool University, Suzhou, China\\
\texttt{Linxiao.Xu@xjtlu.edu.cn}\\[1em]
Preliminary version,
\today
}

\newtheorem{theorem}{Theorem}[section]
\newtheorem{lemma}[theorem]{Lemma}
\newtheorem{proposition}[theorem]{Proposition}
\newtheorem{corollary}[theorem]{Corollary}

\theoremstyle{definition}
\newtheorem{definition}[theorem]{Definition}

\theoremstyle{remark}
\newtheorem{remark}[theorem]{Remark}
\newtheorem{example}[theorem]{Example}

\numberwithin{equation}{section}

\newcommand{\Z}{\mathbb{Z}}
\newcommand{\N}{\mathbb{N}}
\newcommand{\R}{\mathbb{R}}

\newcommand{\SD}{\,\underset{\neq}{\mathrm{dom}}\,}

\newcommand{\sC}{\mathscr{C}}

\DeclareMathOperator{\dom}{\, dom}
\DeclareMathOperator{\PLC}{PLC}
\DeclareMathOperator{\coeff}{coeff}
\DeclareMathOperator{\pos}{Pos}

\DeclareMathOperator{\supp}{supp}
\DeclareMathOperator{\GL}{GL}

\DeclareMathOperator{\rad}{Rad}
\DeclareMathOperator{\cone}{cone}
\DeclareMathOperator{\Int}{Int}
\DeclareMathOperator{\conv}{conv}
\DeclareMathOperator{\aff}{aff}
\DeclareMathOperator{\nonaff}{non-aff}
\DeclareMathOperator{\afftype}{aff-type}


\begin{document}

\begin{abstract}
In this paper we study affine reflection subgroups and their root subsystems in arbitrary infinite Coxeter groups of finite rank.
We give a characterization of limit roots (the accumulation points of the projections of roots onto suitably chosen hyperplanes) associated with affine reflection subgroups.
We describe the relationship between these limit roots and all other limit roots.
In particular, we show that each limit root associated with an affine reflection subgroup may only produce positive values with finitely many positive roots under the bilinear form afforded by the Tits representation (thus satisfying a so called \emph{Pos-finite} condition); and we show that an arbitrary limit root producing positive bilinear form values with only finitely many positive roots must be a convex combination of limit roots associated with affine reflection subgroups. As an application of this Pos-finite condition, we give an independent characterization of when a finitely generated infinite Coxeter group may possess affine reflection subgroups. We also utilize this Pos-finite analysis to partially resolve the question of whether the projections of an infinite injective sequence of roots totally ordered by the dominance partial ordering possesses or not a unique point of accumulation.
\end{abstract}

\maketitle

\section{Introduction}
A \emph{Coxeter system} $(W,R)$ consists of an abstract group $W$ and a generating set $R$ comprised of involutions.
The group $W$, called a \emph{Coxeter group}, is generated by elements of $R$ subject only to braid relations on pairs of generators.
Such a group $W$ can be realized, via the classical \emph{Tits representation}, as a reflection group acting on a real vector space $V$, where the set of \emph{reflections} is taken to mean the set
$T:=\bigcup_{w\in W} w R w^{-1}$.
A \emph{reflection subgroup} of $W$ is a subgroup of $W$ generated by a subset of $T$. It is well known that a reflection subgroup of a Coxeter group is itself a Coxeter group. The classification of the reflection subgroups of Coxeter groups is of considerable importance, and it arises in a multitude of mathematical contexts. Complete classifications of the reflection subgroups of finite Coxeter groups and affine Coxeter groups can be found in \cite{DL11}. However, for general Coxeter groups, the classification of their reflection subgroups is yet to be completed.
Among all reflection subgroups, those isomorphic to affine Coxeter groups are called \emph{affine reflection subgroups}. In this paper we study affine reflection subgroups in Coxeter groups with finitely many generators.

The \emph{root system} of a Coxeter group $W$ is the collection of all \emph{roots} which are representative $-1$-eigenvectors of the reflections in $W$. It is well known that the root system of a Coxeter group $W$ is a discrete set in $V$. However, through an approach  first systematically developed in \cite{HLR11} (whose underlying philosophy may trace back to the work in \cite{VK}), it has been observed that, when dealing with finitely generated infinite Coxeter groups, the projections of their roots onto suitably chosen hyperplanes are always contained in certain compact sets. These projections of roots are called \emph{normalized roots}, and they may exhibit intricate asymptotic behaviours as observed in \cite{HLR11} and \cite{DHR13}. The set of accumulation points of normalized roots for a Coxeter group $W$ is denoted by $E(W)$, and its elements are called the \emph{limit roots} of $W$. It was observed in
\cite{HLR11} that $E(W)$ is contained in the isotropic cone of $V$ with respect to the bilinear form arising from the Tits representation, and there is a natural $W$-action on $E(W)$, and it was later shown in \cite{DHR13} that this action was minimally faithful. In \cite{Dyer12} and \cite{DHR13} it has been further proven that the convex hull of $E(W)$ is the topological closure of the projection of the so-called \emph{imaginary cone} onto the chosen hyperplane in which the normalized roots lie. The notion of an imaginary cone was first introduced in the context of Kac-Moody algebra  as the cone generated by positive imaginary roots (for an accessible reference, please refer to \cite{VK}), and has been later adapted to the setting of Coxeter groups as a (potentially not strict) subset of the so-called \emph{dual of the Tits cone}. An authoritative work on the imaginary cones of Coxeter groups can be found in \cite{Dyer12}. Also in \cite{HE1} one may uncover a treasure trove of related materials.

In this paper, we study the limit roots arising from affine reflection subgroups of a Coxeter group $W$. In the case when $W$ is infinite but finitely generated, we give a characterization of the set of limit roots arising from affine reflection subgroups.
If $\eta$ is one such limit root, then we show that $\eta$ satisfies the following:
\begin{itemize}
\item[(i)] the cardinality of the set $\{x \in \Phi^+\mid (x, \eta)>0 \}$ is finite, where $\Phi^+$ denotes the set of positive roots, and $(\,,\,)$ denotes the bilinear form arising from the Tits representation, henceforth known as the \emph{Pos-finite} condition
(Theorem~\ref{th:fin1});
\item[(ii)] the support of $\eta$  is a connected set and remains connected (in the sense of Coxeter graphs) under the action of $W$ (Proposition~\ref{connect});
\item[(iii)] $\eta$ is in the intersection of the imaginary cone and the normalized isotropic cone (Proposition~\ref{tits}).
\end{itemize}
We prove that if a limit root $\eta$ satisfies (i) and (ii) of the above then $\eta$ must be a limit root arising from an affine reflection subgroup (Theorem~\ref{th:aff2}), and we show that a point in the normalized isotropic cone is a limit root arising from an affine reflection subgroup if and only if it satisfies (ii) and can be acted by an element of $W$ into the fundamental domain of the imaginary cone (Theorem~\ref{rep}). We also show that, if a limit root $\eta$ satisfies (i) but not (ii) then $\eta$ must be a non-trivial linear combination of limit roots arising from affine reflection subgroups. In fact, we give a characterization of when every point in the convex hull of certain limit roots arising from affine reflection subgroups are themselves limit roots (Proposition~\ref{hull}). 

Moreover, by studying the properties of limit roots arising from affine reflection subgroups of finitely generated infinite Coxeter groups, we obtain an independent characterization of when (irreducible) finitely generated infinite Coxeter groups may possess affine reflection subgroups; indeed, in Corollary \ref{co:std}  we observe that irreducible affine reflection subgroups of finitely generated infinite Coxeter groups are precisely the infinite reflection subgroups of conjugates of standard affine parabolic subgroups. Similar observations concerning affine reflection subgroups with three or more generators can be found in \cite[Theorem 3.3]{CM05}.

Furthermore, to the best of our knowledge, the following two natural questions concerning limit roots and the imaginary cone are not explicitly answered in the literature: which points in the isotropic set are necessarily limit roots; and is the containment of the imaginary cone in the dual of the Tits cone proper? It turns out that by studying those limit roots satisfying the Pos-finite condition and those limit roots not satisfying the Pos-finite condition,  we prove that the containment of the imaginary cone in the dual of the Tits cone is strict in the case of finitely generated non-affine infinite Coxeter groups (Corollary~\ref{strict}), with the imaginary cone containing precisely those limit roots satisfying the Pos-finite condition and the complement of the imaginary cone in the dual of the Tits cone containing those limit roots which do not satisfy the Pos-finite condition.  
In this process, we actually obtain a description of which points in the isotropic cone must be limit roots (Proposition~\ref{KQ}).

Finally, we give an application of limit roots corresponding to affine reflection subgroups to partially resolve an open problem of whether or not an infinite sequence of normalized roots whose corresponding positive roots are totally ordered by an important partial ordering, called \emph{dominance}, possesses a unique point of accumulation (conjectured in \cite{DHR13}). The dominance partial ordering, first introduced in \cite{BH93}, was defined on the root systems of infinite Coxeter groups, and has natural interpretation and application in geometric group theory, known as \emph{parallel walls}.
Refinement on existing literature on this partial ordering may supply experts in both representation theory and geometric group theory with further tools to study infinite Coxeter groups. In this paper, we show how the study of limit roots and dominance of finitely generated infinite Coxeter groups may intertwine. Specifically we prove in Proposition~\ref{domlimit} the following result: if any accumulation point of an infinite sequence of normalized roots whose corresponding positive roots are totally ordered by dominance is a limit root arising from affine reflection subgroups, then that sequence of normalized roots is, in fact, convergent, with that affine limit root being the sole point of accumulation.  Moreover, this enables us to construct concrete prototypes of convergent sequences of normalized roots other than the well-known dihedral cases.


The organization of this paper beyond this introduction section is as follows: in Section~$2$, we recall basic definitions and properties of Coxeter groups. It contains a slight generalization to the classical Tits representation of Coxeter groups in a
construction known as a \emph{Coxeter datum} (following the works of \cite{HE2}, \cite{DK94}, \cite{DK09}, \cite{FU1}, and \cite{FU2}).
We recall properties of reflection subgroups in Coxeter groups,
and specifically, we recall the beautiful characterization of \emph{canonical generators} of reflection subgroups, as obtained in
\cite{MD87} and \cite{MD90}.  In Section~$3$, we recall a characterization of the root systems of affine Coxeter groups taken from \cite{VK}. Also we recall an elegant result from \cite{DK94} and \cite{DK09} establishing a Coxeter group to be affine if the bilinear form associated with its Tits representation admits a non-zero radical.
We also collect a few basic properties of  the dominance partial ordering defined in the root systems of infinite Coxeter groups
(as contained in \cite{BH93}, \cite{BB98}, \cite{FU1}, and \cite{FU2}).
In Section~$4$, we recall the construction of normalized roots and limit roots as pioneered in \cite{HLR11}.
In Section~$5$, we present the new results of this paper as described in the earlier part of this introduction.
For readers familiar with the works in \cite{BN68}, \cite{HE2}, \cite{MD87}, \cite{MD90}, \cite{HM}, \cite{BH93}, \cite{DK94}, \cite{DK09}, and \cite{HLR11} they may commence directly from Section~$5$. In section~$6$, we give an application of limit roots relating to affine reflection subgroups to the study of infinite injective sequence of roots totally ordered by dominance. In Proposition~\ref{domlimit} we partially resolve a conjecture raised in \cite{DHR13}. We also give some characterization of when a trio of roots ordered by dominance generate an affine reflection subgroup.

We would like to sincerely thank Professor R.~B.~Howlett, Professor G.~Lehrer, and Professor M.~Dyer for their valuable suggestions, some of these suggestions have directly given rise to numerous results presented in this paper, as well as correcting a number of errors contained in an earlier version of this paper. Finally but not least, we would like to express our sincere gratitude to the referees of this paper for their valuable advice leading to a number of corrections and improvements contained in the current version of this paper.

\section{Definitions and notations}

For a set $A$, define $\PLC(A)$, the \emph{positive linear combination of $A$}, to be the set given by
$$\{\,\sum\limits_{a\in A'} \lambda_a a\mid \text{$A'\subseteq A$, $A'$ is finite and non-empty,
$\lambda_{a}>0$ for all $a\in A'$}\,\};$$
and we define $\cone(A)$, the \emph{positive cone spanned by $A$},  to be the set
$$\cone (A) =\PLC(A)\cup \{0\};$$
and we use the notation of $\conv(A)$ to denote the \emph{convex hull} of $A$.

Suppose that $(W, R)$ is a Coxeter system in the sense of \cite{HH81} or \cite{HM} with
$R$ being a set of involutions generating $W$
subject only to the condition that
$(st)^{m_{st}}=1$ for all $s\neq t\in R$ for some $m_{st}\in \N_{\geq 2}\cup\{\infty\}$.
We stress that the same abstract Coxeter group $W$ may be realized in a number of different ways, that is, different Coxeter systems
may give rise to the same abstract group with distinct generating sets, and possibly more importantly, different sets of reflections. For an elementary example one may consider the dihedral group of order $20$ and the direct product of $\Z_2$ with the dihedral group of order $10$, in which the same abstract group can be realized differently as Coxeter groups with incompatible sets of reflections.  Having this in mind, in this paper, in most cases, we take the preference of referring to a Coxeter system $(W, R)$ instead of just referring to a Coxeter group $W$; and when we simply use the expression ``let $W$ be a Coxeter group'', the reader should interpret it as a statement concerning not only the group $W$, but more importantly the particular presentation specified by a Coxeter system $(W, R)$.

Following the tradition of \cite{DK94} and \cite{DK09}, we make the following definition of a root basis:
\begin{definition}
 \label{def:datum}
Given a Coxeter system $(W, R)$, let $V$ be a vector space over $\R$ and let $(\,,\,)$ be a symmetric bilinear
form on $V$ and let $\Pi=\{\alpha_s\mid s\in R\}$ be a set of vectors of $V$ in bijective correspondence with $R$. Then $\Pi$ is called a \emph{root
basis} if the following hold:
\begin{itemize}
 \item [(C1)] $(\alpha_s, \alpha_s)=1$ for all $s\in R$, and if $s, t$ are distinct generators
then either $(\alpha_s, \alpha_t)=-\cos(\pi/m_{st})$ if $m_{st}$ is a finite integer, or else $(\alpha_s, \alpha_t) \leq -1$ if $m_{st}=\infty$;
 \item [(C2)] $0\notin \PLC(\Pi)$.
 \end{itemize}
\end{definition}

When $\Pi$ is a root basis the triple $\mathscr{C}=(\,V, \, \Pi,
\,(\,,\,)\,)$ is called a \emph{Coxeter datum}.
For each non-isotropic $a\in V$ (that is, $(a, a)\neq 0$), let $\rho_a$ be the \emph{reflection} corresponding to $a$; that is
$\rho_a \in \GL(V)$  is defined by $\rho_a x=x-2\frac{(x, a)}{(a, a)}a$, for all $x\in V$.
%
%
Let $G_{\mathscr{C}}$ be the subgroup of $\GL(V)$ generated by
$\{\,\rho_{\alpha_s}\mid s\in R\,\}=\{\,\rho_a\mid a\in \Pi\,\}$.
%
%
Let $\phi_{\mathscr{C}}\colon W\to G_{\mathscr{C}} $ be the \emph{Tits representation}. It is the unique isomorphism with
$$\phi_{\mathscr{C}}\colon W\to G_{\mathscr{C}}, \text{ satisfying $\phi_{\mathscr{C}}(r_a)=\rho_a$ for all $a\in \Pi$ }.$$
Under the Tits representation there is a  $(\,,\,)$-preserving
$W$-action on $V$: for
each $w\in W$ and $x\in V$, define $wx\in V$ by $wx=\phi_{\mathscr{C}}(w)x$.
We call $(W, R)$ the \emph{abstract Coxeter system} associated to the Coxeter datum $\mathscr{C}$, and
we call $W$ a \emph{Coxeter group} of rank $\#R$ (where $\#$ denotes
cardinality). We use $\ell \colon W\to \N$ to denote the \emph{length function} of $W$ (with respect to $R$).

Let $\Phi\subset V$ be the \emph{root system} of the Coxeter system $(W, R)$ arising form the Coxeter datum $\mathscr{C}$.
That is,
$$\Phi=W\Pi =\{wa\mid w\in W, \, a\in \Pi\}.$$
The set $\Phi^+=\Phi\cap \PLC(\Pi)$ is called the set of \emph{positive roots},
and the set $\Phi^-=-\Phi^+$ is called  the set of \emph{negative roots}.

When $\Pi$ is a basis for the space $V$, for each $x\in V$ there is a unique expression
$$x=\sum_{a\in \Pi}\coeff_a(x) a, \quad \coeff_a(x)\in \R, \forall a\in \Pi,$$
and we define the \emph{support} of $x$ to be
$$\supp(x)=\{a\in \Pi\mid \coeff_a(x)\neq 0\}.$$

Define $T=\bigcup_{w\in W} w Rw^{-1}$. We call $T$ the \emph{set of reflections}
in $W$. For each $x\in \Phi$ we have $x=w \alpha_s$ for
some $w\in W$ and $s\in R$ (or equivalently, $\alpha_s\in \Pi$).
Define $r_x =wsw^{-1}\in T$. Then $\phi_{\mathscr{C}}(r_x)=\rho_x$, and we call $r_x$ the \emph{reflection corresponding to $x$}.
It is readily checked that $T=\{r_x \mid x\in \Phi\}$.


Define functions $N\colon W\to \mathcal{P}(\Phi^+)$ and
$\overline{N}\colon W\to \mathcal{P}(T)$ (where $\mathcal{P}$ denotes power set)
by requiring
$N(w)=\{\,x\in \Phi^+\mid wx\in \Phi^-\,\}$ and
$\overline{N}(w)=\{\,t\in T \mid \ell(wt)<\ell(w)\,\}$ for all $w\in W$.
Note that then $\overline{N}(w)=\{\,r_x\mid x\in N(w)\,\}$.
\begin{definition}
A subgroup $W'$ of $W$ is a \emph{reflection subgroup} of $W$ if  $W'=\langle
T'\rangle$ for some $T'\subseteq T$.
If $W'$ is a reflection subgroup of $W$, we write $W'\leq W$,
furthermore, we let
\begin{align*}
R(W')&=\{\,t\in T\mid \overline{N}(t)\cap W'=\{t\}\,\}\\
\noalign{\hbox{and}}
\Pi(W')&=\{\,x\in \Phi^+\mid r_x\in R(W')\,\}.
\end{align*}
We call $R(W')$ the set of \emph{canonical generators} of $W'$, and we call
$\Pi(W')$ the set of \emph{canonical roots} of $W'$.
\end{definition}
It was shown by Dyer (\cite[Theorems 3.3]{MD90}) and Deodhar (\cite{VD82} and \cite{VD89}) that $(W', R(W'))$
is a Coxeter system whenever $W'$ is a reflection subgroup of $W$, and  furthermore, \cite[Theorem 4.4]{MD90} gave a characterization of
$\Pi(W')$: For $\Delta \subseteq \Phi^+$, let $W'=\langle r_x\mid x\in \Delta\rangle$. Then $\Delta =\Pi(W')$ if and only if
$$
(a, b)\in \{\,-\cos(\pi/n)\mid \text{$n\in \N$ and $n\geq 2$}\,\}\cup (-\infty,
-1],
$$
whenever $a, b\in \Delta$ are distinct. We shall use this characterization quite often in the rest of this paper.


%

Let $(\,,\,)'$ be the restriction of $(\,,\,)$ on the subspace
$\R(\Pi(W'))$ of $V$ spanned by $\Pi(W')$.  Then $\mathscr{C'}=(\,\R(\Pi(W')),\,
\Pi(W'),\, (\,,\,)'\,)$ is a Coxeter datum with $(W', R(W'))$ being the
associated abstract Coxeter system. Thus, the notion of a root system applies to
$\mathscr{C'}$. We let $\Phi(W')$, $\Phi^+(W')$, and $\Phi^-(W')$ be,
respectively, the set of roots, positive roots, and negative roots for the datum
$\mathscr{C'}$. Then $\Phi(W')=W'\Pi(W')$ and \cite[Theorem 3.3]{MD90}
yields that
$$\Phi(W')=\{x\in \Phi\mid r_x \in W'\},$$
and we call $\Phi(W')$ the \emph{root subsystem} of the reflection subgroup $W'$.  We denote
$\Phi^+(W')=\Phi(W')\cap \PLC(\Pi(W'))$ and $\Phi^-(W')=-\Phi^+(W')$,
the sets of positive and negative roots in this subsystem. In this paper a
reflection subgroup $W'$ is called a \emph{dihedral reflection subgroup} if
$\#R(W')=2$. Moreover, we distinguish two types of infinite dihedral reflection subgroups:
\emph{infinite affine reflection subgroups} and \emph{infinite non-affine dihedral reflection subgroups}.
Specifically, if $W'$ is an infinite dihedral reflection subgroup of $W$ and $\Pi(W')=\{a, b\}\subset \Phi^+$ then
$$(a, b)\leq -1.$$
If $(a, b)=-1$, then $W'$ is an infinite affine dihedral reflection subgroup; whereas if $(a, b)<-1$, then $W'$ is an infinite non-affine dihedral reflection subgroup. Note that the two types of infinite dihedral reflection groups are each isomorphic to the abstract infinite dihedral Coxeter group $\widetilde{A}_1$.


A \emph{standard parabolic subgroup} $W_M$ of $W$ is defined as follows, for $M\subseteq \Pi$, let $R_M:=\{r_a\mid a\in M \}$,
and $W_M :=\langle R_M \rangle$.
Observed that $\Phi(W_M)=\Phi\cap \R M$. We may also use the notation $W_{R_M}$ in place of $W_M$. A \emph{parabolic subgroup} of $W$ is any conjugate of a standard
parabolic subgroup of $W$.

\section{Some background results}
\label{sec:basics}


\begin{definition}
Let $(W, R)$ be a Coxeter system with $\mathscr{C} =(V, \Pi, (\,,\,))$ being an associated Coxeter datum.
We define the following two sets related to $(\,,\,)$:
\begin{itemize}
\item[(i)] the \emph{radical} of $(\,,\,)$, denoted by $\rad$ is given by
$$\rad :=\{v\in V\mid (v, u)=0 \text{ for all $u\in V$}\};$$
\item[(ii)] the \emph{isotropic cone} of $(\,,\,)$, denoted by $Q$ is given by
$$Q:=\{v\in V\mid (v, v)=0\}.$$
\end{itemize}
\end{definition}

It is well-known that in an irreducible affine Coxeter group $W$, the bilinear form $(\,,\,)$ afforded by the Tits representation of $W$ has a radical with the following properties:

\begin{proposition}\textup{(\cite[Proposition 2.6\,and\,Theorem 2.7]{HM})}
\label{hum}
Let $(W, R)$ be a Coxeter system in which $W$ is an irreducible affine Coxeter group and $R$ is a finite set of generators,
and let $\mathscr{C}=(V, \Pi, (\,,\,))$ be an associated Coxeter datum.
Then the bilinear form $(\,,\,)$ has a radical $\rad$ such that
\begin{itemize}
\item[(i)] $\dim (\rad) =1$;
\item[(ii)] the radical $\rad$ coincides with the isotropic cone $Q$.\qed
\end{itemize}

\end{proposition}

\begin{lemma}\textup{(\cite[Proposition 6.3]{VK})}
\label{lem: affine}
Let $(W, R)$ be a Coxeter system in which $W$ is an irreducible affine Coxeter group (in the sense of \cite{HM}), with  associated Coxeter datum $\mathscr{C}=(V, \Pi, (\,,\,))$ and  corresponding root system $\Phi$. Then each element of $\Phi$ is congruent modulo $\rad$ to an element
of the root subsystem corresponding to some fixed finite standard parabolic subgroup of $W$.
\qed
\end{lemma}

The following is a key result which affords an elegant method to ascertain whether a Coxeter group $W$ is affine or not, based on the radical associated with its Coxeter datum $\mathscr{C}=(V, \Pi, (\,,\,))$:
\begin{lemma}\textup{(\cite[Lemma 6.1.1]{DK09})}
\label{lem:krammer}
Let $W$ be an irreducible infinite Coxeter group of finite rank. Then $W$ is non-affine implies that
$$\rad\cap \cone(\Pi)=\{0\}.$$
\qed
\end{lemma}

The next result is also well known:
\begin{proposition}\textup{(\cite[Lemma 4.9]{HLR11})}
\label{pp:findi}
Suppose that $W$ is a Coxeter group of finite rank. Then $W$ is infinite if and only if $W$ contains an infinite dihedral reflection subgroup.
\qed
\end{proposition}

\begin{definition}
\label{def:dom} Let $W$ be a Coxeter group and let $W'$ be a reflection subgroup of $W$ with $x, y\in
\Phi(W')$.

\noindent\rm{(i)}\quad We say that $x$ \emph{dominates} $y$ with respect to $W'$ if
$$\{\,w\in W'\mid wx\in \Phi^-(W')\,\}\subseteq \{\,w\in W'\mid wy\in
\Phi^-(W')\,\}.$$ If $x$ dominates $y$ with respect to $W'$ then we write
$x\dom_{W'} y$, furthermore, we write $x\SD_{W'} y$ if $x\dom_{W'} y$ with $x\neq y$, and we simply write $ x\dom y$ in case of $x \dom_W y$.

\noindent\rm{(ii)}\quad Define
$D(x) =\{\, y\in \Phi^+\mid \text{ $x\SD y$\,}\}$.
If $x\in\Phi^+$ and $D(x)=\emptyset$ then we call $x$ \emph{elementary}. For each non-negative integer $n$, define $D_{n}=\{\,x\in \Phi^+\mid
\#D(x)=n\,\}$.
\end{definition}

It is clear from the above definition that for each $x\in \Phi^+$, the set $D(x)$ is finite; indeed, if $x\dom y$ and $y\in \Phi^+$ then $r_x y\in \Phi^-$, and whence $y\in N(r_x)$, a finite set of cardinality $\ell(r_x)$.

It was shown in \cite{BH93} by Brink and Howlett that when a Coxeter group $W$ is finitely generated
then the set of elementary roots is finite. This finiteness property then enabled Brink and Howlett to establish
that all finitely generated Coxeter groups are automatic. Brink later gave an complete construction of $D_0$ for
all such Coxeter groups in \cite{BB98}. Subsequently in \cite{FU1}, it was shown that in a finitely generated Coxeter group $W$
the sets $D_n$ are finite for all $n\in \N$, and furthermore, each $D_n$ is non-empty for all infinite Coxeter groups.


It is readily checked that dominance with respect to any reflection subgroup $W'$ of
a Coxeter group $W$ is a partial ordering on $\Phi(W')$.
The following lemma summarizes some basic properties of
dominance:
\begin{lemma}{\textup{(\cite[Lemma 2.2]{FU1})}}
  \label{lem:basicdom}
\rm{(i)}\quad Let $x, y \in \Phi$ be arbitrary. Then there is dominance between  $x$ and $y$ if and only
if $(x,y) \geq 1$.

\noindent\rlap{\rm{(ii)}}\qquad Dominance is $W$-invariant, that is,  if $x \dom y$ then
 $wx \dom wy$ for all $w \in W$.

\noindent\rlap{\rm{(iii)}}\qquad Let $x, y \in \Phi$ be such that $x \dom y$.
Then
    $-y \dom -x$.
\qed
\end{lemma}

We close this section with the following easy but useful result.
\begin{lemma}\label{lem3.11}
Suppose that $a, b\in \Phi^+$ with $(a, b)\leq -1$, and suppose that $x\in \Phi$ with $x\dom a$.
Then $(x, b)\leq -1$.
\end{lemma}
\begin{proof}
Since $(a, b)\leq -1$, it follows that $a\dom -b$. Thus we have $x\dom a\dom -b$, and hence $(x, b)\leq -1$.
\end{proof}

\section{Limit roots}
Throughout this section we fix a Coxeter system $(W, R)$ and an associated Coxeter datum $\mathscr{C}=(V, \Pi, (\,,\,))$ with $\Phi$ being the corresponding root system. Furthermore, we impose that $\#\Pi=\#R<\infty$. It is well known that $\Phi$ is a discrete set (with respect to the standard topology), since $W$ acts discretely on $V$. Nevertheless, following the approach of studying the so-called
\emph{normalized roots} as pioneered in \cite{HLR11}, we may study the asymptotic behaviours that these normalized roots exhibit if we consider the elements of $\Phi$ as representatives of the directions whose corresponding reflections generate the reflection group $W$.
In \cite{HLR11} the following novel approach has been adopted to obtain such a set of representative directions: consider a projective version of the infinite root system $\Phi$ by cutting the cone $\PLC(\Pi)$, in which the positive roots are located, by a suitable affine hyperplane $\mathscr{H}$; and in doing so, we obtain the so-called \emph{normalized root system} $\widehat{\Phi}$, consisting of the intersections of the rays spanned by the roots in $\Phi^+$ with the hyperplane $\mathscr{H}$. The great advantage of this construction is that the properties of a root system then ensures that these normalized roots are contained in the convex hull of $\widehat{\Pi}$ (consisting of the intersections of the rays $\{\R x\}_{x\in \Pi}$ with the hyperplane $\mathscr{H}$). Under the assumption that $\#\Pi$ being finite, the convex hull of $\widehat{\Pi}$ forms a compact polytope, and thus allowing $\widehat{\Phi}$ to possess a non-empty set of accumulation points, and by studying these accumulation points we may gain a glimpse of
the distribution of roots in infinite root systems. For this construction to work properly the hyperplane $\mathscr{H}$ should conform with the following:


\begin{definition}\textup{(\cite{HLR11}, \cite{DHR13})}
\label{def: trans}
An affine hyperplane $V_1$ of codimension $1$ in $V$ is
called \emph{transverse} to $\Phi^+$  if for each $a\in \Pi$ the ray $\R_{>0}a$  intersects $V_1$ in
exactly one point. Given a hyperplane $V_1$ transverse to
$\Phi^+$, let $V_0$ be the hyperplane that is parallel to $V_1$ and contains the origin.
\end{definition}

\begin{remark}
Recall the standard result that $\Phi=\Phi^+\biguplus\Phi^-$ (where $\biguplus$ denotes
disjoint union). This combined with the requirement that $0\notin \PLC(\Pi)$ immediately imply that it is
always possible to find a hyperplane containing the origin that separates $\Phi^+$ and $\Phi^-$ (see \cite[5.2]{HLR11} for more details).
By suitably translating this hyperplane it is always possible to find a hyperplane transverse to  $\Phi^+$.
\end{remark}

Let $V_1$ be a transverse hyperplane and let $V_0$ be as in the preceding definition.
Let $V_0^{+}$ be the open half space induced by $V_0$ that contains $V_1$. Observe that
$V_0^{+}$ contains  $\PLC(\Pi)$. Since
$\Phi^+\subset \PLC(\Pi)\subset V_0^{+}$, and $V_1$ is parallel to the boundary of
$V_0^{+}$, it follows that
$$\#(V_1\cap \R_{>0} \beta) =1$$ for each $\beta\in \Phi^+$.


\begin{definition}
\label{df: norm}
Let $V_1$ be a transverse hyperplane  in $V$, and let $V_0$ be obtained from $V_1$ as in Definition~\ref{def: trans}.
\begin{enumerate}

\item For each $v\in V\setminus V_0$,  the unique intersection point of $\R v$
and the transverse hyperplane $V_1$  is denoted $\widehat{v}$. The
\emph{normalization map} $\pi_{V_1}\colon V\setminus V_0\to V_1$ is given by $\pi_{V_1}(v):=\widehat{v}$ for all $v\in V\setminus V_0$.
Set $\widehat{\Phi}:=\pi_{V_1}(\Phi)$, and the elements
$\widehat{x}\in\widehat{\Phi}$ are called \emph{normalized roots}.


\item Let $|\cdot|_1\colon V\to \R$ be the unique linear map satisfying the requirement that
$|v|_1=0$ for all $v\in V_0$, and $|v|_1=1$ for all $v\in V_1$.
\end{enumerate}
\end{definition}

Observe that $\pi_{V_1}(-x)= \pi_{V_1}(x)$ for all $x\in V$, and  $\widehat{y} =\frac{y}{|y|_1}$ for all $y\in V\setminus V_0$.

Also observe that $\widehat{\Phi}\subseteq\conv(\widehat{\Pi})$, (recall that
$\conv(X)$ denotes the convex hull of a set $X$), and $\widehat{\Pi}=\{\,\widehat{x}\mid x\in \Pi \,\}$.
Since $W$ is finitely generated, $\Pi$ is a finite set, and then we see that $\widehat{\Phi}$ is contained in the compact set
$\conv(\widehat{\Pi})$. Consequently, if $\Pi$ is finite, then the accumulation points of $\widehat{\Phi}$
are contained in $\conv(\widehat{\Pi})$.

\begin{definition}
Keep all the notations of the previous definition.
\begin{enumerate}
\item The set of \emph{limit roots}  $E(W)$ of the Coxeter group $W$ (with respect to $V_1$)
 is the set of accumulation points of
$\widehat{\Phi}$. For a reflection subgroup $W'\leq W$, the set of limit roots $E(W')$ arising from
the reflection subgroup $W'$ (with respect to $V_1$) is the set of accumulation points of $\widehat{\Phi(W')}$.

\item The \emph{normalized isotropic cone}  $\widehat{Q}$ (associated to the bilinear form $(\,,\,)$ and with respect to $V_1$) is the set $\widehat{Q}= Q\cap V_1$.
 \end{enumerate}
\end{definition}

\begin{example}
\label{eg: special}
Suppose that $\Pi$ forms a basis for the space $V$.
Then the hyperplane
$V_1:=\{\,v\in V \mid \sum_{a\in \Pi} \coeff_a(v)=1\,\}$ is a transverse hyperplane, and
$V_0:=\{\, v\in V \mid \sum_{a\in \Pi} \coeff_a(v)=0    \,\}$ is the corresponding hyperplane obtained from translating $V_1$
to contain the origin. Observe that under these conditions the corresponding $|\cdot|_1$ has the property that
$|v|_1=\sum_{a\in \Pi} \coeff_a(v)$ for all $v\in V$.
\end{example}


\begin{remark}
\label{rmk:std}
For the rest of this section and the whole of the next section, we shall adopt the set up as in the preceding example, that is, we shall assume that $\Pi$ forms a basis for $V$, and we shall take
$$V_1=\{\,v\in V\mid \sum_{a\in \Pi} \coeff_a(v)=1\,\}.$$
\end{remark}

Next we have an elegant result from \cite{HLR11}:
\begin{theorem}{\textup{(\cite[Theorem 2.7\,and\,Proposition 2.14]{HLR11})}}
\label{thm: limit}
Let $V_1$ be a transverse hyperplane, and
let $\widehat{Q}$ be the associated normalized isotropic cone. Then
\begin{enumerate}
\item $E(W)\neq \emptyset$ if and only if the  Coxeter
group $W$ is infinite.

\item $E(W)\subseteq \widehat{Q}$.\qed
\end{enumerate}
\end{theorem}
%



Following the convention set in \cite{HLR11}, we define
$$D=\bigcap_{w\in W} w(V\setminus V_0)\cap V_1 = V_1\setminus \bigcup_{w\in W} w V_0,$$
and we define the $\cdot$-action of $W$ on $D$ as follows: for any $w\in W$ and $x\in D$,
$$w\cdot x =\widehat{wx}.$$
Observe that the property that $WD=\bigcup_{w\in W}w D \subseteq V\setminus V_0$
guarantees that this $\cdot$-action of $W$ on $D$ is well-defined. Furthermore,
we observe that each $w\in W$ acts continuously on $D$. The next result is taken from \cite{HLR11}
which summarizes a number of key facts:

\begin{proposition}\textup{(\cite[Proposition 3.1]{HLR11})}
\label{pp:D}
Let $V_1$ be a transverse hyperplane, and let $\widehat{\Phi}$ and $E(W)$ be the corresponding normalized roots and limit roots.
\begin{itemize}
\item[(i)] $\widehat{\Phi}$ and $E(W)$ are contained in $D$.
\item[(ii)] $\widehat{\Phi}$ and $E(W)$ are stable under the $\cdot$-action of $W$;
moreover
$$\widehat{\Phi}= W\cdot \Pi.$$
\item[(iii)] The topological closure $\widehat{\Phi}\uplus E(W)$ is stable under
the $\cdot$-action of $W$.
\end{itemize}
\qed
\end{proposition}

Furthermore, it has been observed in \cite{HLR11} that the $\cdot$-action has the following nice geometric description:

\begin{proposition}\textup{(\cite[Proposition 3.5]{HLR11})}
\label{pp:geom}
Keep previous notations.
\begin{itemize}
 \item[(i)] Let $\alpha\in \Phi$, and $x\in D\cap Q$. Denote by $L(\widehat{\alpha}, x)$ the line containing
           $\widehat{\alpha}$ and $x$. Then
					\begin{itemize}
					\item[(a)] if $(\alpha, x)=0$, then $L(\widehat{\alpha}, x)$ intersects $Q$ only at $x$, and $r_{\alpha}\cdot x =x$;
					\item[(b)] if $(\alpha, x)\neq 0$, then $L(\widehat{\alpha}, x)$ intersects $Q$ in two distinct points, namely, $x$ and
					           $r_{\alpha}\cdot x$.
					\end{itemize}
\item[(ii)] Let $\alpha_1$ and $\alpha_2$ be two distinct roots in $\Phi$, $x\in L(\widehat{\alpha_1}, \widehat{\alpha_2})\cap Q$, and $w\in W$. Then $w\cdot x\in L(w\cdot \alpha_1, w\cdot \alpha_2)\cap Q$.					
\end{itemize}
\qed

\end{proposition}

Given a transverse hyperplane $V_1$, the limit roots coming from a given rank $2$ reflection subgroup can
be observed inside $E(\Phi)$. Take two distinct positive roots $a$ and $b$, denote by $W'=\langle r_a, r_b\rangle$
the dihedral reflection subgroup of $W$ generated by the reflections $r_a$ and $r_b$ corresponding to the two positive roots.
Let $\Pi(W')=\{a', b'\}$ be the set of canonical roots for the dihedral reflection subgroup $W'$. Then \cite[Theorem 4.4]{MD90} yields that
$(a', b')\in (-\infty, -1]\cup \{-\cos(\pi/n)\mid n\in \N_{\geq 2}\}$. It is clear that $0\notin \PLC(\{a', b'\})$, since $a', b'\in \PLC(\Pi)$
and $0\notin\PLC(\Pi)$ by the definition of $\Pi$. Now let $V'=\R a'\oplus\R b'$, a two dimensional subspace of $V$ spanned by $\Pi'=\{a', b'\}$, and let $(\,,\,)'$ be the restriction of $(\,,\,)$ on $V'$. Then $\mathscr{C}'=(V', \Pi', (\,,\,)')$ is also a Coxeter datum with associated root system
$\Phi' = W' \Pi'$. Observe that the hyperplane $V_1$ is also transverse with respect to $\Phi'$. Let $E(W')$ denote the set of limit roots of $W'$ in the root system $\Phi'$ with respect to $V_1$. Then the following was observed in $2.3$ of \cite{HLR11}:

\begin{proposition}{\textup{(\cite{HLR11})}}
\label{pp:dih}
Given the set up in the preceding paragraph,
\begin{itemize}
\item[(i)]$E(W')=Q\cap L(\widehat{a'}, \widehat{b'})=E(W)\cap L(\widehat{a'}, \widehat{b'})$;
\item[(ii)]the cardinality of $E(W')$ is $0$, $1$, or $2$, respectively, precisely when
$|(a', b')|<1$, $|(a', b')|=1$, or $|(a', b')|>1$.\qed
\end{itemize}
\end{proposition}

It turned out by considering the cardinality of $E(W)$ we may easily ascertain whether $W$ is affine or not.

\begin{proposition}\textup{(\cite[Corollary 2.16]{HLR11})}
\label{pp:singleton}
Suppose that $W$ is an irreducible Coxeter group. Then $E(W)$ is a singleton set if and only if $W$ is an affine Coxeter group.
\qed
\end{proposition}

The last two propositions allow us to deduce that every infinite dihedral reflection subgroup of an irreducible affine Coxeter group must be affine. Conversely, as a corollary to \cite[Theorem 6.2 (a)]{DHR13}, one can establish that an infinite non-affine Coxeter group must possess an infinite dihedral reflection subgroup that is non-affine. Thus, we have:
\begin{proposition}
\label{prop:aff}
Let $W$ be an irreducible Coxeter group of finite rank. Then $W$ is affine if and only
if every infinite dihedral reflection subgroup of $W$ is affine.
\qed
\end{proposition}

To close this section on necessary introductory material on limit roots, we further recall the following important result which establishes that when $W$ is a finitely generated irreducible and infinite Coxeter group,  the limit set $E(W)$ itself is a $W$-invariant object (with respect to the $\cdot$-action of $W$), and this action is minimally faithful:

\begin{theorem}\textup{(\cite[Proposition 3.2]{HLR11}) \cite[Theorem 3.1]{DHR13}}
\label{th:limit}
Suppose that $(W, R)$ is an irreducible Coxeter system in which $|R|<\infty$ and $|W|=\infty$. Then
\begin{itemize}
\item[(i)] $|wx|_1>0$ for all $w\in W$, and $x\in E(W)$;
\item[(ii)] If $x\in E(W)$, then $E(W)=\overline{W\cdot x}=\overline{\{\,w\cdot x\mid w\in W\,\}}$, that is, the $\cdot$-action of $W$ on $E(W)$ is minimal.\qed
\end{itemize}
\end{theorem}
With more technical but otherwise predictable generalization, it can be observed that the irreducibility requirement in part (i) of Theorem~\ref{th:limit} can be removed.

\section{Limit roots relating to affine reflection subgroups}

Throughout this section, each Coxeter system $(W, R)$ is understood to have an associated Coxeter datum
 $\sC=(V, \Pi, (\,,\,))$ with $\Pi$ being a finite set, and $\Phi$ being the corresponding root system.
Further, unless otherwise stated the Coxeter group $W$ is assumed to be infinite.
Keeping the set-up as in Remark~\ref{rmk:std}, specifically, we assume that $\Pi$ forms a basis for $V$, and we let $V_1$ be the transverse hyperplane as in Remark~\ref{rmk:std}, and for any reflection subgroup $W'\leq W$, let $E(W')$ be the set of limit roots of $W'$ with respect to the transverse hyperplane $V_1$.


A number of results in this section, including Theorem~\ref{th:fin1}, Proposition~\ref{fundd}, and Proposition~\ref{rep}, can be seen as consequences of results contained in \cite{Dyer12}. More specifically, in \cite{Dyer12}, a novel approach of realizing the Tits cone of a Coxeter group $W$ in the same vector space $V$ bearing the Tits representation of $W$ (rather than in the algebraic dual of $V$, as in the more conventional treatment) has been adopted; and a certain functoriality condition was established on the imaginary cones of a Coxeter group $W$ and the reflection subgroups of $W$, namely, the imaginary cone of a reflection subgroup is contained in the imaginary cone of the ambient group $W$. These results are significant and they often require working in set-ups different from the more conventional ones, and therefore, for the purposes of self-containedness and ease of reading, we have included independent proofs using more standard methodologies. We sincerely thank an anonymous referee who pointed out these alternative possibilities to us.

Suppose that $W'$ is an irreducible affine reflection subgroup of $W$, and let $\Pi(W') =\{\, a_1, a_2, \ldots, a_s\,\}$.
Let $\eta$ be the unique limit root in $E(W')$.
By Lemma~\ref{lem: affine}, the limit root $\eta$ is in the radical of the bilinear form $(\,,\,)$ restricted to the subspace
$V_{W'}$ spanned by $a_1, a_2, \ldots, a_s$,
and it then follows that $(\eta, x)=0$ for all $x\in \Phi(W')\subset V_{W'}$.
However, it is possible that $(\eta, x)>0$ for roots $x\in \Phi^+\setminus \Phi(W')$, and such positive roots may potentially
form an infinite set. In this section, amongst other things, we shall prove that $(\eta, x)>0$ for only finitely many positive roots $x\in \Phi^+$. Note that Propositions~\ref{pp:findi},~\ref{pp:singleton}, and~\ref{prop:aff} together imply that if
$W'$ is an irreducible affine reflection subgroup of $W$ then $W'$ has an affine dihedral reflection subgroup $W''\subseteq W'$ with
$E(W'')=E(W')$. Hence many discussions on general affine reflection subgroups, especially their limit roots, can be simplified to discussions on affine dihedral reflection subgroups.

Before we give a characterization of those limit roots in $E(W)$ arising from affine reflection subgroups, let us first look at some limit roots that possibly do not arise from affine reflection subgroups. The first candidate for such possibly non-affine limit roots might come from infinite non-affine dihedral reflection subgroups.

Choose $a, b\in \Phi^+$ such that $(a, b)=-\cosh\theta<-1$ for some $\theta >0$, and let $W'=\langle r_a, r_b\rangle$ be the dihedral reflection subgroup generated by the reflections corresponding to $a$ and $b$. Note first that $W'$ is infinite and non-affine. Proposition~\ref{pp:dih} then yields that $E(W')$ consists of two distinct limit roots. Direct calculations then show that the
intersection of the  isotropic cone $Q$ with the subspace
$\R a\oplus \R b$ consists of two lines
\begin{align*}
&\R((\cosh\theta+\sinh\theta)a+b)\\
\noalign{\hbox{and}}
&\R((\cosh\theta-\sinh\theta)a+b)
        = \R(a+(\cosh\theta+\sinh\theta)b),
\end{align*}				
and if we let $\eta_1$ and $\eta_2$ be as following, then $E(W') =\{\eta_1, \eta_2\}$:
\begin{align*}
\eta_1:=&\frac{(\cosh\theta+\sinh\theta)|a|_1}{(\cosh\theta +\sinh\theta)|a|_1+|b|_1}\widehat{a}+
                 \frac{|b|_1}{(\cosh\theta +\sinh\theta)|a|_1+|b|_1}\widehat{b},\\
								\noalign{\hbox{and}}\\
\eta_2:=&\frac{(\cosh\theta-\sinh\theta)|a|_1}{(\cosh\theta -\sinh\theta)|a|_1+|b|_1}\widehat{a}+
                 \frac{|b|_1}{(\cosh\theta -\sinh\theta)|a|_1+|b|_1}\widehat{b}.
								\end{align*}									
Furthermore,
\begin{align*}
\eta_1 &= \lim_{i\to\infty}\widehat{(r_a r_b)^i a}\in \R((\cosh\theta+\sinh\theta)a+b)\\
\noalign{\hbox{and}}
\eta_2 &= \lim_{i\to \infty}\widehat{(r_b r_a)^i b} \in \R((\cosh\theta-\sinh\theta)a+b).
\end{align*}
And it follows readily that
\begin{align*}
(\eta_1, a) &=\frac{\sinh\theta}{(\cosh\theta +\sinh\theta)|a|_1+|b|_1} >0;\\
\noalign{\hbox{and}}
(\eta_2, b) &=\frac{\sinh\theta (\cosh\theta-\sinh\theta)}{(\cosh\theta -\sinh\theta)|a|_1+|b|_1}>0.
\end{align*}
Note that $(r_b r_a)^i\cdot \eta_1 =\eta_1$ and $(r_a r_b)^i \cdot \eta_2 =\eta_2$ for all $i\in \N$.
Consequently, for all $i\in \N$
\begin{equation}
\label{eq:inf}
(\eta_1,\, (r_a r_b)^i a)>0 \text{ and }
(\eta_2,\, (r_b r_a)^i b)>0.
\end{equation}
In particular,
\begin{equation*}
\text{$(\eta_1, x)>0$ and
$(\eta_2, y)>0$ for infinitely many $x, y\in \Phi^+$. }
\end{equation*}

However, the above discussion does not immediately rule out the possibility that an affine dihedral reflection subgroup and
an infinite non-affine dihedral reflection subgroup share the same limit root, a situation
as illustrated in the following schematic diagram.
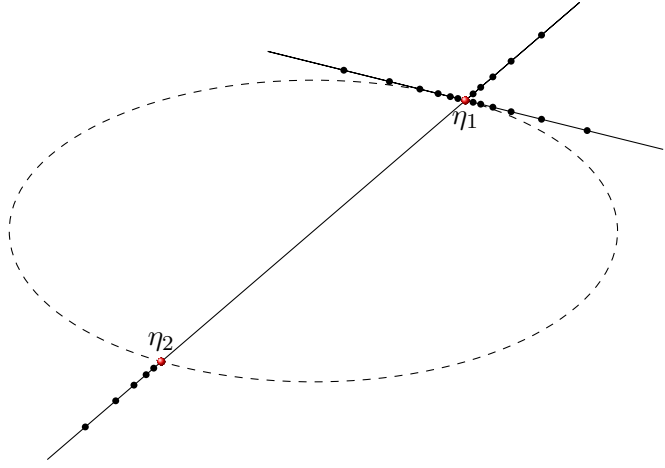
\begin{center}
\begin{tikzpicture}[scale=2]
  \draw[dashed] (0,0) ellipse (2 and 1);
  \coordinate (A) at (1, 0.8660254038);
  \coordinate (B) at (-1, -0.8660254038);
  \draw (A)--+(0.75, 0.8660254038*0.75)--(B)--+(-0.75,-0.8660254038*0.75);
  \draw (A)--+(-1.3, 0.25*1.3)--+(1.3, -0.25*1.3);
  \shade[ball color=red](A) circle (0.025) node[below] {\small{$\eta_1$}};
  \shade[ball color=red](B) circle (0.025) node[above] {\small{$\eta_2$}};
  \fill (1.05, 0.8660254038*1.05) circle (0.65pt);
  \fill (1.1, 0.8660254038*1.1) circle (0.65pt);
  \fill (1.18, 0.8660254038*1.18) circle (0.65pt);
  \fill (1.3, 0.8660254038*1.3) circle (0.65pt);
  \fill (1.5, 0.8660254038*1.5) circle (0.65pt);

  \fill (-1.05, -0.8660254038*1.05) circle (0.65pt);
  \fill (-1.1, -0.8660254038*1.1) circle (0.65pt);
  \fill (-1.18, -0.8660254038*1.18) circle (0.65pt);
  \fill (-1.3, -0.8660254038*1.3) circle (0.65pt);
  \fill (-1.5, -0.8660254038*1.5) circle (0.65pt);
  \fill (1, 0.8660254038)+(0.05, -0.05*0.25) circle (0.65pt);
  \fill (1, 0.8660254038)+(0.1, -0.1*0.25) circle (0.65pt);
  \fill (1, 0.8660254038)+(0.18, -0.18*0.25) circle (0.65pt);
  \fill (1, 0.8660254038)+(0.3, -0.3*0.25) circle (0.65pt);
  \fill (1, 0.8660254038)+(0.5, -0.5*0.25) circle (0.65pt);
  \fill (1, 0.8660254038)+(0.8, -0.8*0.25) circle (0.65pt);

  \fill (1, 0.8660254038)+(-0.05, 0.05*0.25) circle (0.65pt);
  \fill (1, 0.8660254038)+(-0.1, 0.1*0.25) circle (0.65pt);
  \fill (1, 0.8660254038)+(-0.18, 0.18*0.25) circle (0.65pt);
  \fill (1, 0.8660254038)+(-0.3, 0.3*0.25) circle (0.65pt);
  \fill (1, 0.8660254038)+(-0.5, 0.5*0.25) circle (0.65pt);
  \fill (1, 0.8660254038)+(-0.8, 0.8*0.25) circle (0.65pt);
\end{tikzpicture}
\end{center}
In this schematic  diagram the dotted circle represents the normalized isotropic cone,  and the normalized root subsystems of
two infinite dihedral reflection subgroups, one affine and the other non-affine, are contained in the two straight lines,
with the black dots schematically representing normalized roots, and the two red dots representing possible limit roots.

It turned out that this situation will not arise. Indeed, in this section we will give a characterization of the set of
limit roots arising from affine reflection subgroups. For each $\eta$ arising from an affine reflection subgroup of an infinite Coxeter group of finite rank, we establish where in $\widehat{Q}$ such an $\eta$ must have been located; and we establish the finiteness of the  cardinality of the set
$\{\,x\in \Phi^+\mid (x, \eta)>0\,\}$;
we also establish that the support of $\eta$ is connected and remains connected under the $\cdot$-action of $W$ (recall that a set $M\subseteq \Pi$ is \emph{disconnected} if there are non-empty proper subsets $M_1$ and $M_2$ of $M$ with $M=M_1\uplus M_2$
such that $(a, b) =0$ for each $a\in M_1$ and $b\in M_2$, otherwise $M$ is said to be \emph{connected} and in which case we may say that $M$ forms a connected component of $\Pi$);
moreover, we establish that a limit root simultaneously satisfy these finiteness and connectedness properties must be arising from an affine
reflection subgroup and no other limit roots may satisfy simultaneously these two properties.

\begin{definition}
Let $(W, R)$ be a Coxeter system in which $W$ is an infinite Coxeter group.
A limit root $\eta\in E(W)$ is called an \emph{affine limit root} (or is simply called \emph{affine}) if
there exists an irreducible affine reflection subgroup $W'\le W$ with $E(W')=\{\eta\}$, and we use the set
$$E_{\aff}:=\{\,\eta\in E(W')\mid \text{$W' \leq W$ is affine}\,\}$$
to denote the set of all affine limit roots in $W$.
Also we use the following set to denote the set of all non-affine dihedral limit roots in $W$:
$$E_{\nonaff}^2:=\{\,\eta\in E(W')\mid \text{ $W' \leq W$ is non-affine dihedral} \,\}.$$
\end{definition}
Note that the calculations immediately before the schematic diagram show that if $\eta\in E_{\nonaff}^2$ then $(\eta, x)>0$ for infinitely many $x\in \Phi^+$,
and consequently if we could show that $(\eta, x)>0$ only for finitely many $x\in \Phi^+$ whenever $\eta\in E_{\aff}$, then we would have established that $E_{\aff}\cap E_{\nonaff}^2=\emptyset$, and thus ruling out the possibility raised in the diagram.
To achieve this, we need some preparatory work first.

\begin{definition}
\label{k}
For each $v\in V $, define $\pos(v)\subseteq \Phi^+$ by
$$\pos(v)=\{x\in \Phi^+\mid (v, x)>0\}.$$
Further, define
\begin{align*}
\mathscr{K}&:=\{v\in\cone(\Pi)\mid \pos(v)=\emptyset\}\\
           &\,\,=\{v\in\cone(\Pi) \mid (v, a)\leq 0, \text{ for all $a\in \Phi^+$}\}.
\end{align*}
\end{definition}

Next, we recall the following fundamental result for infinite Coxeter groups of finite rank which is going to be a key ingredient in our proof that $E_{\aff}\cap E_{\nonaff}^2=\emptyset$.

\begin{lemma}\textup{(\cite[Proposition 4.5.5]{ABFB})} 
\label{pp:lower}
The  set
$$\{\,(\alpha, \beta)\mid \alpha, \beta\in \Phi, \,\,        |(\alpha, \beta)|<1 \,\}$$
is finite.
In particular, there exists a fixed $\epsilon >0$ such that $|(\alpha, \beta)|>\epsilon$ whenever $\alpha, \beta \in \Phi$ satisfy $(\alpha, \beta)\neq 0$.
\qed
\end{lemma}

\begin{remark}
Note that by Propositions~\ref{pp:findi}, ~\ref{pp:singleton}, and~\ref{prop:aff},  $\eta\in E_{\aff}$ if and only if there exists an irreducible affine dihedral reflection subgroup  $W'\le W$ with $E(W')=\{\eta\}$.
\end{remark}

Now we are ready to prove that $(\eta, x)>0$ only for finitely many $x\in \Phi^+$ whenever $\eta\in E_{\aff}$:

\begin{theorem}
\label{th:fin1}
If $\eta\in E_{\aff}$, then $\pos(\eta)$ is a finite set.

\end{theorem}
\begin{proof}
Let $D\le W$ be an affine dihedral reflection subgroup with $E(D)=\{\eta\}$ and let $a,b\in\Phi^+$ be such that
$\Pi(D)=\{a,b\}$.
 Then  $(a, b)=-1$
and $(\eta, a)=(\eta, b)=0$.
Letting $a_k=k(a+b)+a$ and $b_k=k(a+b)+b$ (for $k\in \N_{\geq 0}$), we have that
$\widehat{a_k}\to \eta$ and $\widehat{b_k}\to\eta$, and
the root subsystem of $D$ satisfies
$\Phi^+(D)=\{\,a_k, b_k \mid k\ge 0\}$.
 Note also that $a_k\dom a$ and $b_k\dom b$  for all $k\ge 0$.

 Let $x\in\pos(\eta)$.
 Since $\widehat{a_k}\to \eta$ and $\widehat{b_k}\to\eta$, for sufficiently large $k$, we have
 $(x,a_k)>0$ and $(x,b_k)>0$. Fixing such a $k$, we have
 \begin{align*}
 0<(x,a_k)+(x,b_k)=(x,(2k+1)(a+b))=(2k+1)(x,a+b).
 \end{align*}
 Therefore, $(x,a+b)>0$.
 Observe that
  $(x,a)>-1$, since otherwise by Lemma~\ref{lem:basicdom}, $x \dom -a \dom -a_k$ for all $k\in \N$, which implies that $(x, a_k)\leq -1$ for all $k\in \N$,
	but this means that $(x, \eta)\leq 0$, contradicting $x\in \pos(\eta)$. Similarly, $(x, b)> -1$.
	

 We now observe that  there exists $\epsilon>0$ such that
 $(x,a+b)>\epsilon$ for all $x\in\pos(\eta)$. To see this, let $F$ denote the finite set $$F=\{\,(\alpha, \beta)\mid \alpha, \beta\in \Phi, \,\,        |(\alpha, \beta)|<1 \,\},$$
 and define $\epsilon'=\min\{r+s\mid r,s\in F, r+s>0 \}$ and $\epsilon''=1+\min F$.
 Setting $\epsilon=\min(\epsilon',\epsilon'')>0$, we have
 $$
 (x,a+b)=(x,a)+(x,b)\geq \epsilon>0.
 $$
Therefore, since
\begin{align*}
(x,a_k)&=k(x,a+b)+(x,a),\\
\noalign{\hbox{and}}
(x,b_k)&=k(x,a+b)+(x,b),
\end{align*}
we may deduce that there is a fixed $M\in\N$ such that
\begin{equation}\label{eqn:xak}
(x,a_k)\ge 1\text{ and } (x,b_k)\ge 1, \text{ $\forall x\in \pos(\eta)$ and $\forall k\ge M$.}
 \end{equation}
Now, for a contradiction, suppose that $\pos(\eta)$ is an infinite set.
With $M$ as above and for some fixed $k\ge M$,  (\ref{eqn:xak}) then implies that there are dominance between the infinitely
many $x\in \pos(\eta)$ and each of $a_k$ and $b_k$. Since a given root can dominate only finitely many positive roots (see the remark immediately after Definition~\ref{def:dom}), it follows that there are infinitely many positive roots $x\in\pos(\eta)$ satisfying
$x\dom a_k$ and $x\dom b_k$, but then each of such $x$ would simultaneously dominate both $a$ and $b$,  and we have a contradiction to Lemma \ref{lem3.11}.



\end{proof}

It turns out that for a limit root $\eta\in E(W)$, the
finiteness of $\pos(\eta)$ is not a sufficient condition to ensure that $\eta \in E_{\aff}$. To ensure that $\eta$ is the limit root of an affine reflection subgroup of $W$, more conditions on $\eta$ are required. The following is an example demonstrating  this
general fact. We gratefully acknowledge the comments made by Professor M.~Dyer on an earlier version of this paper, and
in particular, we thank him for directing our attention to this example, which first appeared in \cite[Example 5.8]{HLR11}.

\begin{example}
\label{eg5.7}
Let $(W, R)$ be a rank $5$ Coxeter system with an associated Coxeter datum $\mathscr{C}=(V, \Pi, (\,,\,))$
such that $\Pi:=\{a, b, c, d, e\}$ is a basis for $V$, and $(\,,\,)$ satisfies the requirement that the only non-zero
values between a pair of simple roots are the following
$$(a, b) =(d, e)=-1, \quad and \quad (b, c)=(c, d)=-\frac{1}{2}.$$
Then a direct calculation shows that for each $n\in \N$,
$$(r_a r_b r_e r_d)^n c= c+ n^2b +n^2 d+(n^2+n)a +(n^2+n)e.$$
Thus $$\eta_0:=\lim_{n\to \infty} \widehat{(r_a r_b r_e r_d)^n c}=\frac{a+b+d+e}{4}\in E(W).$$
Further, note that $\pos(\eta_0)=\emptyset$, but $\eta_0$ is not a limit arising from an affine reflection subgroup of $W$.
Instead, $\eta_0$ can be seen as the \emph{normalized sum} of the limit roots from two mutually orthogonal affine reflection subgroups, namely,  $W_1:=\langle r_a, r_b\rangle$ and $W_2:=\langle r_d, r_e\rangle$.
Moreover, it can be checked that every point on the line segment in the interval between the two affine limit roots $\frac{a+b}{2}$ and
$\frac{d+e}{2}$ is in fact a limit root. For more details on this interesting behaviour, please see \cite[Example 7.12]{Dyer12}.
It should also be pointed out that even though $\supp(\eta_0)=\{\,a, b, d,e\,\}$, it happens that $\eta_0\notin E(\langle r_a, r_b, r_d, r_e\rangle$).

\end{example}

In addition to the finiteness condition
of the respective $\pos$'s,  a complete characterization of affine limit roots also rests on the number of connected components
of their supports. To prepare us for such a characterization, we still need a few observations. We begin by recalling a few results on the imaginary cones associated to Coxeter groups.

\begin{definition}
\label{ima}
Let $(W, R)$ be an arbitrary Coxeter system, and let $\mathscr{C}=(V, \Pi, (\,,\,))$ be an associated
Coxeter datum with $\Phi$ being the corresponding root system. We define the \emph{imaginary cone}, denoted
by $\mathscr{Z}$, to be
$$\mathscr{Z}:=\bigcup_{w\in W} w\mathscr{K}.$$
where we recall $\mathscr{K}=\{v\in \cone(\Pi)\mid \pos(v)=\emptyset\}$, as introduced in Definition~\ref{k}.
\end{definition}

The concept of the imaginary cone was first introduced in \cite{VK} in the context of Kac-Moody Lie algebras
as the pointed cone spanned by the positive imaginary roots. This concept is later generalized to Coxeter groups.
A definitive reference on the imaginary cones of Coxeter group can be found in \cite{Dyer12}.
We will show, amongst other things, that in a finitely generated Coxeter group $W$, the only limit roots in
$\mathscr{Z}$ are precisely those limit roots relating to  affine reflection subgroups of $W$ (more precisely, linear combinations of affine limit roots). Furthermore, we show that every point in the intersection of the imaginary cone and the normalized isotropic cone is, in fact, a limit root (expressible as a linear combination of affine limit roots).

The following characterization of the imaginary cone of a finitely generated Coxeter group is taken from \cite{FU2}:
\begin{proposition}\textup{(\cite[Lemma 4.4 and Proposition 4.22]{FU2})}
\label{pp:im}
Let $\mathscr{Z}$ be the imaginary cone of a finitely generated Coxeter group $W$. Then
$$\mathscr{Z}=\{v\in U^*\mid \text{ $(v, a)\leq 0$ for all but finitely many $a\in \Phi^+$}\}, $$
where $U^*$ is the dual of the Tits cone given by
$$U^*=\bigcap_{w\in W} w (\cone(\Pi)).$$
\qed
\end{proposition}

\begin{remark}
Theorem \ref{th:fin1} above may be alternatively established by considering the imaginary cones associated with the affine reflection subgroups
of a finitely generated Coxeter group. Suppose that $W$ is a finitely generated infinite Coxeter group possessing an affine reflection subgroup $W'$. Then Propositions~\ref{pp:findi},~\ref{pp:singleton}, and~\ref{prop:aff} yield that there exists some affine dihedral reflection subgroup $W''$ of $W'$ (and hence also a reflection subgroup of $W$) with $E(W'')=E(W')=\{\eta\}$. Let $V''$ be the $2$-dimensional subspace spanned by the canonical roots $\Pi(W'')$ of the affine dihedral reflection subgroup $W''$ and let $(\,,\,)''$ be the restriction of the bilinear form $(\,,\,)$ to $V''$. Let $\Phi(W'')$ be the root subsystem of the affine dihedral reflection subgroup $W''$. Let
$$\mathscr{K}_{W''}=\{v\in \cone(\Pi(W''))\mid (v, a)\leq 0, \text{ for all $a\in \Phi^+(W'')$}\},$$
and let $\mathscr{Z}_{W''}=\bigcup_{w\in W''} w\mathscr{K}_{W''}$ be the imaginary cone of $W''$ with respect to the Coxeter datum
$(V'', \Pi(W''), (\,\,)'')$. Then $$\{\eta\}=\mathscr{K}_{W''}\subset \mathscr{Z}_{W''}.$$ It is a remarkable result of \cite[Theorem 6.3]{Dyer12} that
$$\mathscr{Z}_{W''}\subseteq \mathscr{Z}.$$
Consequently, $\eta\in \mathscr{Z}$, and it follows from Proposition~\ref{pp:im} immediately that $\#\pos(\eta)<\infty$.
\end{remark}

To establish the connection between limit roots of finitely generated infinite Coxeter groups and the associated imaginary cones, we recall the following key result from \cite{DHR13}:

\begin{theorem}\textup{(\cite[Theorem 2.3]{DHR13})}
\label{EZ}
Setting $Z:=\mathscr{Z}\cap V_1$, then the convex hull of $E(W)$ is the topological closure of $Z$.
\qed
\end{theorem}

From Theorem~\ref{EZ}, it is immediately clear that $E(W)\subseteq \overline{Z}$, where $\overline{\,\cdot\,}$ denotes the topological closure of $\cdot$.

It has been shown in \cite[Lemma 2.4]{DHR13} that if $(W, R)$ is an irreducible Coxeter system
in which $W$ is a non-affine infinite Coxeter group then $\Int(\mathscr{K})\neq \emptyset$ (where $\Int(\,\cdot\,)$ denotes the topological interior of $\cdot$).
Note that if $x\in \Int(\mathscr{K})$ then it is clear from the definition of $\mathscr{K}$ that $(x, a)<0$
for all $a\in \Pi$, and in particular, $(x, x)<0$. Thus, no limit root can be in the interior of $\mathscr{K}$, and indeed,
no limit root can be in the interior of $\mathscr{Z}$.



By combining Proposition~\ref{pp:im}, Theorem~\ref{th:limit}, and Theorem~\ref{th:fin1}, we can establish that affine limit roots
are in the imaginary cone (in fact, every limit root $\eta$ with $\#\pos(\eta)<\infty$ is in the imaginary cone).
\begin{proposition}
\label{tits}
If $(W, R)$ is an irreducible Coxeter system then $E(W)\subseteq U^*$. Moreover, if $\eta\in E(W)$ is such that
$\#\pos(\eta)<\infty$, then $\eta\in \mathscr{Z}$, and in particular, $E_{\aff}\subseteq Z$.
\qed
\end{proposition}

We are now able to establish that the containment of the imaginary cone in the dual of Tits cone is strict in the case of
non-affine infinite Coxeter groups of finite rank.
\begin{corollary}
\label{strict}
Let $(W, R)$ be a Coxeter system in which $W$ is an irreducible infinite non-affine Coxeter group and $R$ is a finite generating set.
Then the imaginary cone $\mathscr{Z}$ is a proper subset of the dual of the Tits cone $U^*$. In particular,
$$ E_{\nonaff}^2 \subseteq U^*\setminus \mathscr{Z}\neq \emptyset.$$
\end{corollary}
\begin{proof}
First, we recall Proposition~\ref{pp:findi} which stated that $W$ is infinite if and only if it contains an infinite dihedral reflection subgroup.
Since $W$ is non-affine, it follows from Proposition~\ref{prop:aff}
that there exists an infinite non-affine dihedral reflection subgroup $W'\leq W$, and so
$\emptyset \neq E(W')\subset E_{\nonaff}^2$.
Now let $\eta\in E_{\nonaff}^2$ be arbitrary. Then the calculations at the beginning of
this section yields that $\#\pos(\eta)=\infty$.
Consequently it follows from Proposition~\ref{pp:im} that $\eta\notin \mathscr{Z}$, and so
$E_{\nonaff}^2 \subseteq U^*\setminus \mathscr{Z}\neq \emptyset$.
\end{proof}

For any finitely generated infinite Coxeter group $W$, it was proven in \cite{FR1} that there is a
$W$-equivariant map from the Davis complex of $W$ into the imaginary cone of $W$, and this map is a homeomorphism
onto its image. The Davis complex has many known topological properties, and the existence of such a map may assist
in the discovery of comparable properties in the imaginary cone, whose general topological features are yet to be explored.
The construction of this map built on earlier works in \cite{DK94} and \cite{DK09} where a similar map from the Davis complex of $W$ to the dual of the Tits cone of $W$ was constructed. It was not entirely clear that the map obtained in \cite{FR1} was in fact an improvement to the map constructed in \cite{DK94} and \cite{DK09}. Indeed, it was not clear that the containment of the imaginary cone in the dual of the Tits cone
was strict. This issue is now resolved by Corollary~\ref{strict} in the case of infinite non-affine Coxeter groups.

Next, we strive to complete a characterization of the limit roots in $E_{\aff}$. So far
we have affirmed that
$$E_{\aff}\subseteq \{\,\eta\in E(W)\mid \#\pos(\eta)<\infty\,\}.$$
We shall show that if the group $W$ is not an affine Coxeter group, then this containment
is strict. Further, we will complete a characterization of the set $\{\,\eta\in E(W)\mid \#\pos(\eta)<\infty\,\}$
itself, and this set is actually closely related to affine limit roots.
Furthermore, we shall establish the connection of this set of limit roots with the imaginary cone $\mathscr{Z}$.
In particular, we are able to ascertain that each point in $\mathscr{Z}\cap \widehat{Q}$ is in fact a limit root.
To complete these tasks, we still need a few observations on limit roots.

\begin{lemma}
\label{lem:posfin}
Let $\eta\in E(W)$ be such that $\#\pos(\eta) <\infty$. Then there exists some
$w\in W$ with $w\eta\in \mathscr{K}$.
\end{lemma}
\begin{proof}
By Propositions~\ref{pp:im} and ~\ref{tits}, the fact that $\#\pos(\eta)<\infty$ implies that $\eta\in Z= \mathscr{Z}\cap V_1$.
Then the definition of the imaginary cone gives that
$\eta \in \bigcup_{w\in W}w \mathscr{K}$. Consequently, there exists some $w\in W$ with
$w\eta \in \mathscr{K}$.
\end{proof}

It turned out that in a Coxeter group $W$ of finite rank, if a limit root $\eta\in E(W)$ satisfies the condition that
$\#\pos(\eta)<\infty$, then the $W$-orbit of $\eta$ contains a unique element in $\mathscr{K}$. This is a consequence
of a more general fact that $\mathscr{K}$ is a \emph{fundamental domain} for the action of $W$ on $\mathscr{Z}$.
That is, the $W$-orbit of every element of $\mathscr{Z}$  contains one and only one element in $\mathscr{K}$.
\begin{proposition}
\label{fundd}
Keep the notations of Definition~\ref{ima}. Let $\eta\in \mathscr{K}$, and let $M:=\{\,a\in \Pi\mid (a, \eta)=0\,\}$.
Then
\begin{equation}
\label{stabilizer}
\{\,w\in W\mid w\eta\in \mathscr{K}\,\}=W_M  =\{w\in W\mid w\eta=\eta\}.
\end{equation}
That is, $\mathscr{K}$ is a \emph{fundamental domain} for the action of $W$ on $\mathscr{Z}$.

In particular, in a Coxeter group $W$ of finite rank, if a limit root $\eta\in E(W)$ satisfies the condition that
$\#\pos(\eta)<\infty$, then the $W$-orbit of $\eta$ contains a unique element in $\mathscr{K}$.
\end{proposition}
\begin{proof}
Note that by the definition of $\mathscr{Z}$, it is clear that the $W$-orbit of every element of $\mathscr{Z}$  contains some
element(s) of $\mathscr{K}$, and we only need to prove (\ref{stabilizer}) to establish the uniqueness.

Clearly $W_M\subseteq \{\,w\in W\mid w\eta=\eta\,\}\subseteq \{\,w\in W\mid w\eta\in \mathscr{K}\,\}$.
Hence it is enough to show that
$\{\,w\in W\mid w\eta\in \mathscr{K}\,\}\subseteq W_M$.

Let $w\in W$ such that $w\eta =\eta'\in \mathscr{K}$. If $w\neq 1$ then we may choose $a\in \Pi$
with $w=w'r_a$ and $\ell(w)=\ell(w')+1$. Note that $w'a\in \Phi^+$, and hence
$$
0\geq (w'a, \eta')= (w'a, w'r_a \eta)
                  = (a, r_a\eta)
									= (r_a a, \eta)
									=-(a, \eta)
									\geq 0,
$$
since $\eta,\,\eta'\in \mathscr{K}$ and $(\,,\,)$ is $W$-invariant.
Consequently $a\in M$.
Setting $w_0=w$ and $w_1 =w'$,  and we have
$\eta' =w\eta= w_0\eta =w'r_a\eta=w_1\eta$,
and $w_0 W_M =w_1 W_M$.
Now if $w_1\neq 1$ then we may repeat the above argument and find some $w_2\in W$ satisfying the following three conditions:
$\ell(w_2)=\ell(w_1)-1=\ell(w_0)-2$; $w_2 \eta =w_1 \eta =w_0\eta=\eta'$; and
$w_2 W_M= w_1 W_M=w_0 W_M$. If $w_2\neq 1$, we can repeat the argument again,
and so construct a sequence $w_0, w_1, w_2, \ldots \in W$. Since the length decrease at
each step this sequence must terminate; however, this process can be continued so long
$w_i\neq 1$. Thus $w_i=1$ for some finite $i$, and
$$wW_M =w_0 W_M =w_1 W_M =w_2 W_M =\cdots =W_M,$$
proving $w\in W_M$.

Finally, note that under the given conditions, $\#\pos(\eta)<\infty$ implies that $\eta\in \mathscr{Z}$, and the rest of the
proposition follows from the above proof.
\end{proof}

\begin{remark}
For those readers who are familiar with the treatment of the Tits cone being contained in the same vector space $V$ bearing the Tits representation as pioneered in \cite{Dyer12}, Proposition~\ref{fundd} above translates into the statement that $\mathscr{K}$ is contained in the negative of the fundamental chamber on $V$ in the set-ups of \cite{Dyer12}.
\end{remark}

%


\begin{lemma}
\label{lem:rad}
Let $x\in \mathscr{K}\cap Q$, and $x\neq 0$.
Let $M:=\supp(x)\subseteq \Pi$.
Then the restriction of the bilinear form $(\,,\,)$ on the subspace spanned by $M$ has a non-zero radical
$\rad(M)$, and $x\in \rad(M)$.
\end{lemma}
\begin{proof}
Write $x=\sum_{a\in M}\coeff_a(x) a$.
Since $0\neq x\in \mathscr{K}\subseteq \PLC(\Pi)$, it follows that $\coeff_a(x)>0$ for all $a\in M$.
Now $x\in Q$ implies that
\begin{equation}
\label{KQ1}
0=(x, x)=\sum_{a\in M}\coeff_a(x)(x, a).
\end{equation}
Since $x\in \mathscr{K}$, it follows that $(x, a)\leq 0$ for all $a\in \Pi$, and hence (\ref{KQ1}) yields that
$(x, a)=0$ for all $a\in M$. Consequently, the restriction of the bilinear form $(\,,\,)$ on the subspace spanned by $M$
has a non-zero radical, $\rad(M)$, and moreover, $x\in \rad(M)$.
\end{proof}

We can deduce that the converse of Lemma \ref{lem:rad} is also true.

\begin{lemma}
\label{lem:rad2}
Let $M\subseteq \Pi$ be such that the restriction of the bilinear form
$(\,,\,)$ on the subspace spanned by $M$ has a non-zero radical $\rad(M)$, and suppose that moreover, $\rad(M)\cap \PLC(M)\neq\emptyset$. Then
$\rad(M)\cap \PLC(M) \subseteq \mathscr{K}\cap Q$.
\end{lemma}
\begin{proof}
Let $x\in \rad(M) \cap \PLC(M)$ be arbitrary. Then $(x, a)=0$ for all $a\in M$, and moreover,
$(x, a')\leq 0$ for all $a'\in \Pi\setminus M$ by the definition of $\Pi$ and the fact that $x\in \PLC(M)$.
Therefore $x\in \mathscr{K}$.

Next, to establish that $x\in Q$, it is enough to observe that  the fact $x\in \rad(M)$ implies
that $(x, x)=\sum_{a\in M}\coeff_a(x)(x,a )=0$.
\end{proof}

\begin{proposition}
\label{pp:rad}
Suppose that $\eta\in E(W)$ such that  $\pos(\eta)$ is a finite set. Then
there exists $M\subseteq \Pi$ such that the restriction of the bilinear
form $(\,,\,)$ on the subspace spanned by $M$ has a non-zero radical $\rad(M)$ satisfying the condition that
$\rad(M)\cap \PLC(M)\neq \emptyset$, and
there exists some $w\in W$ such that
$$w\cdot \eta \in \rad(M).$$
\end{proposition}
\begin{proof}
By Lemma \ref{lem:posfin}, there exists some $w\in W$ such that $w\cdot \eta\in \mathscr{K}$. Note that
then $0\neq w\cdot \eta \in \mathscr{K}\cap Q$. Consequently, it follows from Lemma~\ref{lem:rad}
that there exists some $M\subseteq \Pi$ such that the restriction of the bilinear form $(\,,\,)$ to the subspace spanned
by $M$ has  a non-zero radical $\rad(M)$, with $\rad(M)\cap\PLC(M)\neq \emptyset$, and $w\cdot \eta \in \rad(M)$.
\end{proof}

The next technical result follows from \cite[8.4]{Dyer12}, but for completeness sake we include a proof here.

\begin{proposition}
\label{KQ}
Let $\eta\in \mathscr{K}\cap \widehat{Q}$ be such that $M:=\supp(\eta)$ is connected.
Then $\eta\in E_{\aff}$, and the standard
parabolic subgroup $W_M$ is an irreducible affine reflection subgroup of $W$  with
$E(W_M)=\{\eta\} $.
\end{proposition}
\begin{proof}
Lemma~\ref{lem:rad} yields that the restriction of the bilinear form $(\,,\,)$ on the subspace spanned by $M$ has
a non-zero radical $\rad(M)$ with $\eta \in \rad(M)\cap \PLC(M)$.
Now by Lemma~\ref{lem:krammer}, the conditions $M$ being connected (and obviously finite) and
$\rad(M)\cap \PLC(M)\neq \emptyset$ together imply that the reflection subgroup
$W_M:=\langle r_a\mid a\in M\rangle$ is an irreducible affine Coxeter group. Then Proposition~\ref{hum}
yields that $\rad(M)$ is one dimensional, and it coincides with the isotropic cone restricted to the subspace spanned by $M$.
Consequently, the normalized isotropic cone restricted to the subspace spanned by $M$ consists of a single point, and hence $E (W_M)=\{\eta\}$.

\end{proof}

It turned out that $\mathscr{K}\cap \widehat{Q} \subset E(W)$, and indeed, $\mathscr{Z} \cap \widehat{Q} \subset E(W)$. We shall show these once we have a better understanding of
those limit roots whose supports are not connected.

\begin{lemma}
\label{dim}
Let $M$ be a connected subset of $\Pi$ such that the restriction of the bilinear form $(\,,\,)$
on the subspace spanned by $M$ has a non-zero radical $\rad(M)$ satisfying $\rad(M)\cap \PLC(M)\neq \emptyset$. Then $\dim(\rad(M))=1$.
\end{lemma}
\begin{proof}
Note that $\#M<\infty$, since $\#\Pi<\infty$. Then by Lemma~\ref{lem:krammer}, the conditions $M$ being connected
and $\rad(M)\cap \PLC(M)\neq \emptyset$ establish that the standard parabolic subgroup
$W_M:=\langle r_a\mid a\in M\rangle$ is an irreducible affine Coxeter group. Then Proposition~\ref{hum}
yields that $\rad(M)$ is one dimensional.
\end{proof}

Lemma~\ref{dim} has the following easily seen consequence.

\begin{corollary}
\label{components}
Let $M$ be a finite subset of $\Pi$ such that the restriction of the bilinear form $(\,,\,)$
on the subspace spanned by $M$ has a non-zero radical $\rad(M)$, and there exists some
$x\in \rad(M)\cap\PLC(M)$ with $\supp(x)=M$.
Let $M =M_1 \uplus \cdots \uplus M_m$ be a decomposition of $M$ into connected components.
Then for each $i=1, \ldots, m$, the restriction of the bilinear form $(\,,\,)$
on the subspace spanned by $M_i$ has a $1$-dimensional radical $\rad(M_i)$ with $\rad(M_i)\cap \PLC(M_i)\neq \emptyset$.
Furthermore, $\rad(M) = \rad(M_1)\oplus \cdots \oplus \rad(M_m)$.

\end{corollary}
\begin{proof}
Under the given conditions, we may write $x=x_1+\cdots + x_m$, where $x_i\in \PLC(M_i)$ and $\supp(x_i)=M_i$ for $i=1, \ldots, m$.
Since the $M_i$'s are pairwise disconnected, it follows from the properties of a root basis that $(a, b)=0$  whenever $a\in \R M_i$ and $b\in \R M_j$ with $i\neq j$.
Similarly, if $a\in \Pi\setminus M_i$ then $(x_i, a) \leq 0$;
on the other hand, Lemma~\ref{lem:rad2} ensures that $x\in \mathscr{K}$, and so if $a\in M_i$ then $(x_i, a)=(x, a) \leq 0$. Thus $x_i\in \mathscr{K}$ for each $i$.
Now since $0=(x, x) =\sum_{i=1}^m (x_i, x_i)$, and $(x_i, x_i)= \sum_{a\in M_i}|\coeff_a(x_i)| (x_i, a)\leq 0$,
it follows that $x_i \in Q$ for each $i$.
Then Lemma~\ref{lem:rad} yields that for each $i$, the restriction of the bilinear form $(\,,\,)$
on the subspace spanned by $M_i$ has a non-zero radical $\rad(M_i)$ such that $0\neq x_i \in\rad(M_i)\cap \PLC(M_i)$. Because each $M_i$ is connected, Lemma~\ref{dim} yields that $\dim(\rad(M_i))=1$.
Finally,  it is readily checked that $\rad(M) = \rad(M_1)\oplus \cdots \oplus \rad(M_m)$.

\end{proof}

\begin{proposition}
\label{connect}
Suppose that $\eta\in E(W)$ is a limit root arising from an irreducible affine reflection subgroup $W'\leq W$. Then
$\supp(\eta)$ is a connected subset of $\Pi$.
\end{proposition}
\begin{proof}
Note that Proposition~\ref{pp:singleton} and Proposition~\ref{prop:aff} yield that we may choose an affine dihedral
reflection subgroup $W''\leq W'$ such that $\{\eta\}=E(W'')$. Set $\Pi(W'')=\{a, b\}$. Then $a, b\in \Phi^+$ and
$(a, b)=-1$.
Then $\supp(\eta) =\supp(a)\cup \supp(b)$. Since $a$ and $b$ are roots, it can be easily checked that $\supp(a)$ and
$\supp(b)$ are both connected. Observe then that the fact $(a, b)=-1$ implies that $\supp(a)\cup \supp(b)$ is also connected.
For otherwise,
\begin{align*}
a&=\sum_{\alpha\in \supp(a)}\coeff_{\alpha}(a) \alpha, \text{ where all the $\coeff_{\alpha}(a) >0$},\\
\noalign{\hbox{and}}
b&=\sum_{\beta\in \supp(b)}\coeff_{\beta}(b) \beta, \text{ where all the $\coeff_{\beta}(b) >0$},
\end{align*}
and if $\supp(a)\cup\supp(b)$ is not connected, then note that $(\alpha, \beta)=0$ for all $\alpha\in \supp(a)$ and for all
$\beta\in \supp(b)$. But then
$$(a, b)=\sum_{\alpha\in \supp(a)} \sum_{\beta\in \supp(b)}\coeff_{\alpha}(a) \coeff_{\beta}(b) (\alpha, \beta)=0,$$
contradicting that $(a, b)=-1$.
\end{proof}

\begin{theorem}
\label{th:aff2}
Let $\eta\in E(W)$ be such that $\#\pos(\eta)<\infty$ and $\supp(w\cdot\eta)$ is connected for each $w\in W$. Then there exists an irreducible affine reflection subgroup $W'$ of
$W$ satisfying $E(W')=\{\eta\}$.
\end{theorem}
\begin{proof}
By Lemma~\ref{lem:posfin}, the fact  $\#\pos(\eta)<\infty$   implies that there exists $w\in W$ such that
$$0\neq w\cdot \eta\in \mathscr{K}\cap Q.$$
Let $M:=\supp(w\cdot\eta)$, and note that $M$ is connected.
Then Proposition~\ref{KQ} yields that $W_M:=\langle r_a\mid a\in M\rangle\leq W$
is an irreducible affine standard parabolic subgroup, and $\{w\cdot \eta\}=E(W_M)$.
Then it follows that $w^{-1}W_M w$, being a conjugate of an irreducible affine reflection subgroup of $W$,
is itself an irreducible affine reflection subgroup of $W$, and $\{\eta\}=E(w^{-1}W_M w)$.

\end{proof}

Combining Theorem~\ref{th:fin1}, Proposition~\ref{connect}, and Theorem~\ref{th:aff2} we immediately have the following characterization of affine limit roots in
an infinite Coxeter group of finite rank.

\begin{theorem}
\label{th:key}
Let $\eta\in E(W)$. Then
$\eta\in E(W')$ for some affine reflection subgroup $W'$ of $W$ if and only if $\pos(\eta)$ is a finite set and
$\supp(w\cdot\eta)$ is connected for each $w\in W$.
\qed
\end{theorem}

In fact, the set of affine limit roots also admits the following alternative characterization:
\begin{theorem}
\label{rep}
If $\eta\in E(W)$, then
$\eta\in E_{\aff}$ if and only if there exists
some $w\in W$ such that $w\cdot \eta \in \mathscr{K}$ with $\supp(w\cdot \eta)$ being connected.
\end{theorem}
\begin{proof}
Note that under the given conditions, if $\eta\in E_{\aff}$ then Theorem~\ref{th:fin1}
yields that $\#\pos(\eta)<\infty$, and Proposition~\ref{lem:posfin} then yields that there exists some $w\in W$ with
$w\cdot \eta\in \mathscr{K}$, and moreover, Proposition~\ref{connect} ensures that $\supp(w\cdot\eta)$ is connected.

Conversely, let $\eta\in E(W)$ be such that there exists some $w\in W$ with $w\cdot \eta\in \mathscr{K}\cap Q$ and $M:=\supp(w\cdot \eta)$ being connected. Then Proposition~\ref{KQ} yields that $W_M$ is an irreducible affine reflection subgroup of $W$ and $\{w\cdot\eta\}=E(W_M)$.
Consequently, $w^{-1}W_M w$ is an irreducible affine reflection subgroup of $W$ with $\{\eta\} =E(w^{-1}W_M w)\subseteq E_{\aff}$.
\end{proof}

\begin{remark}
For those readers who are familiar with the treatment of realizing the Tits cone in the same vector space bearing the Tits representation, Theorem~\ref{rep} above can be proved alternatively: If $\eta\in E$ such that $\pos(\eta)$ is finite, then there exists $w\in W$ such that $w\eta$ is in the negative of the Tits cone realized in $V$. Furthermore, since $w\cdot\eta\in E\subset \PLC(\Pi)$, it follows that $w\eta\in \PLC(\Pi)\cap \mathscr{K}$. Then \cite[8.4 and 8.8]{Dyer12} readily establishes that $w\eta$ is a convex combination of limit roots from $E_{\aff}$, and consequently, so is $\eta$.
\end{remark}

\begin{corollary}
\label{co:std}
An infinite Coxeter group $W$ of finite rank contains an affine reflection subgroup if and only if $W$ contains an affine standard parabolic subgroup. In particular,
irreducible affine reflection subgroups in $W$ are precisely  the irreducible infinite reflection subgroups of  affine parabolic subgroups of $W$.
\end{corollary}
\begin{proof}
Suppose that $W$ has an affine standard parabolic subgroup. Then there is nothing to prove, for a standard parabolic subgroup of $W$ is, \emph{a priori}, a reflection subgroup of $W$.

Conversely, suppose that $W$ has an irreducible affine reflection subgroup $W'$, and suppose that $\{\eta\}=E(W')$. Then
Theorem~\ref{rep} yields that there exists some $w\in W$ such that $w\cdot \eta \in \mathscr{K}$ and $N:=\supp(w\cdot\eta)$
is connected. It follows from Proposition~\ref{KQ} that $W_N$ is an irreducible affine reflection subgroup of $W$.
Next, let
$$M:=\{\,a\in \Pi\mid (a, w\cdot\eta)=0\,\}.$$
Note that Proposition~\ref{fundd} then yields that $N\subseteq M$.
Now, if $a\in M\setminus N$, then $(a, b)=0$ for all $b\in N$, since otherwise
there exists $b\in N$ with $(a, b)<0$, and since
$$w\cdot \eta =\coeff_b(w\cdot \eta) b +\sum_{c\in N\setminus\{b\}}\coeff_c(w\cdot \eta) c,$$
with all the coefficients being strictly positive,
then it follows from the definition of a root basis that
\begin{align*}
(a, w\cdot \eta) &=(a, \coeff_b(w\cdot \eta) b +\sum_{c\in N\setminus\{b\}}\coeff_c(w\cdot\eta) c\,)\\
                 &=\coeff_b(w\cdot\eta) (a, b) + \sum_{c\in N\setminus\{b\}} \coeff_c(w\cdot\eta) (a, c)\\
								 &<0,
\end{align*}
contradicting that $a\in M$. Thus $M=N\uplus N'$ where $N$ and $N'$ are disconnected subsets of $\Pi$,
and $W_M$ is the direct product of standard parabolic subgroups $W_N$ and $W_{N'}$.
Now for any $w'\in W'$, we have
\begin{align*}
(ww'w^{-1})\cdot (w\cdot\eta)=(w w' w^{-1} w)\cdot \eta=(w w')\cdot\eta &=w\cdot(w'\cdot\eta)\\
                                                                        &=w\cdot \eta.
\end{align*}
That is, $w W' w^{-1}$ stabilizes $w\cdot\eta\in \mathscr{K}$. Then it follows from Proposition~\ref{fundd} that
$w W' w^{-1}\subseteq W_M$, and indeed, $w W' w^{-1}\subseteq W_N$.
\end{proof}
Observations comparable to Corollary~\ref{co:std} may be found in \cite[Corollary 8.8]{Dyer12} and \cite[Theorem 3.3]{CM05} (ultimately tracing back to \cite{DK94}, \cite{DK09}, and other unpublished works of Daan Krammer) for affine reflection subgroups of ranks greater than or equal to three. Corollary~\ref{co:std} provides a unified approach to deal with arbitrary affine reflection subgroups.

\begin{proposition}
\label{sum}
Suppose that $\eta\in E(W)$ such that $\#\pos(\eta)<\infty$ and $\eta\notin E_{\aff}$. Then $\eta$ is in the convex hull
of limit roots from $E_{\aff}$.
\end{proposition}
\begin{proof}
Since $\#\pos(\eta)<\infty$, it follows from Lemma~\ref{lem:posfin} that there exists some $w\in W$ such that
$w\cdot \eta\in \mathscr{K}\cap Q$.
Let $M:=\supp(w\cdot \eta)$. Then Lemma~\ref{lem:rad} gives that the bilinear form $(\,,\,)$ restricted to the subspace spanned by $M$ has a non-zero
radical $\rad(M)$ such that $w\cdot \eta\in \rad(M)\cap \PLC(M)$.
Since $\eta\notin E_{\aff}$, it follows that $w\cdot\eta\notin E_{\aff}$, and then
Proposition~\ref{fundd} and Theorem~\ref{rep} together yield that $M$ is not connected.
Now suppose that $M=M_1\uplus \dots \uplus M_m$, where $M_1, \ldots, M_m$ form a complete list of connected components of $M$.
Then Corollary~\ref{components} and Lemma~\ref{lem:krammer} together imply that the standard parabolic subgroups
$W_{M_1}, \ldots, W_{M_m}$ are all irreducible affine reflection subgroups.
Let $\{\eta_i\}:=E(W_{M_i})\subset E_{\aff}$ for each $i=1, \ldots, m$.
Then it follows from Corollary~\ref{components} that $w\cdot \eta$ can be expressed as
a positive linear combination of these $\eta_i$'s.
That is, $w\cdot\eta =\sum_{i=1}^m \lambda_i \eta_i$ for some $\lambda_i \geq 0$ with $\sum_{i=1}^m \lambda_i=1$.
Consequently,
$$\eta =\sum_{i=1}^m \lambda'_i w^{-1}\cdot\eta_i, \text{~where~} \lambda'_i={|w\eta|}_1{|w^{-1}\eta_i|}_1\lambda_i,$$
where $w^{-1}\cdot \eta_i$ is the limit root of the affine reflection subgroup $w^{-1} W_{M_i}w$ for each
$1, \ldots, m$.

\end{proof}

The following is a well-known result, and for completeness we  include a proof here.
\begin{lemma}
\label{supp}
Let $W$ be a finitely generated irreducible Coxeter group which may be either finite or infinite.
Then there exists a positive root $x$ such that $\supp(x)=\Pi$.
\end{lemma}
\begin{proof}
We give an inductive algorithm to compute one such positive root $x$.

Since $W$ is irreducible, its Coxeter graph is connected. We start with $x_1=\alpha$ where
$\alpha$ is an arbitrarily chosen simple root.
Suppose that we already have $x_k\in \Phi^+$ such that $\#\supp(x_k) =k <\# \Pi$.
Since $W$ is irreducible, it follows that there exists
some $\beta\in \Pi\setminus \supp(x_k)$ and $\gamma \in \supp(x_k)$ with $(\beta, \gamma)<0$.
Then set $x_{k+1}:=r_{\beta} x_k \in \Phi^+$ and $\#\supp(x_{k+1}) =k+1$.
We may continue this process until we have found an $x_m$ with $\supp(x_m) =\Pi$.
\end{proof}

\begin{proposition}
\label{hull}
Let $(W, R)$ be a Coxeter system in which $W$ is an irreducible  infinite non-affine Coxeter group and $R$ is a finite
generating set. Let $M_1, \ldots, M_m$ be connected but pairwise disconnected subsets of $\Pi$  which further satisfy the requirement that  their corresponding standard parabolic subgroups $W_{M_1}, \ldots, W_{M_m}$ are irreducible affine reflection subgroups.
Let $\eta_1, \ldots, \eta_m$ be the limit roots of $W_{M_1}, \ldots, W_{M_m}$ respectively. Then
$$\{\,\sum_{i=1}^m \lambda_i \eta_i \mid \text{$\sum_{i=1}^m \lambda_i =1$ and $0\leq \lambda_i \leq 1$ for all $i=1, \ldots, m$ }\,\}\subseteq E(W).$$
That is, each point in the convex hull of the affine limit roots $\eta_1,\ldots, \eta_m$ is a limit root.
\end{proposition}
\begin{proof}
Fix a positive root $x$ with $\supp(x) = \Pi$ whose existence is guaranteed by Lemma~\ref{supp}.

For each $i=1, \ldots, m$, let $\Phi(W_{M_i})$ denote the root subsystem of the affine reflection subgroup $W_{M_i}$.
For each $M_i$, fix a simple root $\alpha_i\in M_i$. By \cite[Proposition 3.13]{HT97}, we can find some $\beta_i \in \Phi^+(W_{M_i})$
satisfying $\beta_i \SD\alpha_i$. Since $W_{M_i}$ is affine, it follows that $(\alpha_i, \beta_i) =1$.
Let $\gamma_i :=r_{\alpha_i} \beta_i\in \Phi^+(W_{M_i})$, and
$$(\alpha_i, \gamma_i)=(\alpha_i, r_{\alpha_i} \beta_i)=-(\alpha_i, \beta_i)=-1.$$
Then by \cite[Theoerm4.4]{MD90}, $\{\alpha_i, \gamma_i\}$ is the set of canonical generators for the dihedral reflection subgroup
$\langle r_{\alpha_i}, r_{\gamma_i}\rangle$. Set $\delta_i =\alpha_i +\gamma_i \in \PLC(\Pi)$.
Note that $(\delta_i, \alpha_i)=(\delta_i, \gamma_i) =0$.
Direct calculations show that
$$\Phi^+(\langle r_{\alpha_i}, r_{\gamma_i}\rangle) =\{\,n\alpha_i +(n\pm 1)\gamma_i\mid n\in \N\,\},$$
and hence $\alpha_i +n\delta_i \in \Phi^+(\langle r_{\alpha_i}, r_{\gamma_i}\rangle)$ for all $n\in \N$.
Thus $\widehat{\delta_i}$ is the sole limit root in $E(\langle r_{\alpha_i}, r_{\gamma_i}\rangle)$.
Since $W_{M_i}$ is an irreducible affine reflection subgroup and $\langle r_{\alpha_i}, r_{\gamma_i}\rangle \subseteq W_{M_i}$,
it follows that $\eta_i =\widehat{\delta_i}$, and in particular, $(\delta_i, a)=0$ for all $a\in M_i$.
Next, note that for each $\ell\in \N$
\begin{align*}
r_{\alpha_i+\ell\delta_i} (x)&=x-2(x, \alpha_i+\ell\delta_i) (\alpha_i+\ell\delta_i)\\
                          &=x-2(x, \alpha_i)\alpha_i -2\ell(x,\delta_i)\alpha_i-2\ell (x, \alpha_i) \delta_i -2\ell^2(x, \delta_i)\delta_i.
\end{align*}
Observe that this is a positive root (since $\supp(x)=\Pi$). Now we let $p_i:=-2(x, \alpha_i)$ and $q_i:=-2(x, \delta_i)$.
Note that $q_i >0$ (since $W$ is irreducible).
Thus the above can be written as
\begin{equation}
\label{af}
r_{\alpha_i+\ell\delta_i} (x) =x+ p_i \alpha_i +\ell q_i \alpha_i +\ell p_i \delta_i +\ell^2 q_i \delta_i.
\end{equation}
Since  $M_1, \ldots, M_m$ are pairwise disconnected subsets of $\Pi$, it follows that $(x_i, x_j)=0$ for all $x_i\in \R M_i$ and $x_j\in \R M_j$, whenever $i\neq j$. Thus from (\ref{af}) we see that applying the pairwise commuting reflections corresponding to the positive roots
$\alpha_1+\ell_1 \delta_1, \,\ldots,\, \alpha_m+\ell_m\delta_m$ to $x$ (where $\ell_1, \ldots, \ell_m\in \N$) produces the following positive root
\begin{align*}
r_{\alpha_m+\ell_m\delta_m} \cdots r_{\alpha_i+\ell_i\delta_i}\cdots r_{\alpha_1+\ell_1\delta_1} (x)=x &+\sum_{i=1}^m p_i \alpha_i
                                                                                              +\sum_{i=1}^m \ell_i q_i \alpha_i\\
																																															&+\sum_{i=1}^m \ell_i p_i \delta_i
																																															 +\sum_{i=1}^m \ell_i^2 q_i \delta_i.
\end{align*}
It can be checked that for any $t_i\in (0, 1)$ with $\sum_{i=1}^m t_i=1$, we may find integral sequences $(\ell_{1_n}), \ldots, (\ell_{m_n})$
with
$$\lim_{n\to \infty} \ell_{i_{n}} =\infty \text{ and }\lim_{n\to \infty}\frac{\ell_{i_n}^2 }{\sum_{i=1}^m \ell_{i_n}^2 }=t_i.$$ Consequently,
$$\{\,\sum_{i=1}^m \lambda_i \eta_i \mid \text{$\sum_{i=1}^m \lambda_i =1$ and $0\leq \lambda_i \leq 1$ for all $i=1, \ldots, m$ }\,\}\subseteq E(W).$$
\end{proof}

\begin{theorem}
\label{ZQE}
Let $W$ be an irreducible Coxeter group.
Then
$$\mathscr{Z} \cap\widehat{Q} \subseteq E(W).$$
\end{theorem}
\begin{proof}
Let $\eta\in \mathscr{Z}\cap \widehat{Q}$ be arbitrary. Then by the definition of the imaginary cone, there
exists some $w\in W$ such that $w \cdot \eta \in \mathscr{K}\cap \widehat{Q}$.
Let $M:=\supp(w\cdot \eta)$. Then Lemma~\ref{lem:rad} gives that the bilinear form $(\,,\,)$ restricted to the
subspace spanned by $M$ has a non-zero radical $\rad(M)$ with $w\cdot \eta\in \rad(M)$.

If $M$ is connected, then Proposition~\ref{KQ} yields that
$w\cdot \eta\in E_{\aff}$, and this in turn means that $\eta\in E_{\aff}$.

If $M$ is disconnected, then let $M=M_1\uplus \cdots\uplus M_m$ be a decomposition of $M$ into
connected components. Then Corollary~\ref{components} and Lemma~\ref{lem:krammer} together establish
that the standard parabolic subgroups $W_{M_1}, \ldots, W_{M_m}$ are all irreducible affine reflection subgroups of $W$.
Denote the limit roots of $W_{M_1}, \ldots, W_{M_m}$ by $\eta_1, \ldots, \eta_m$ respectively.  Then Corollary~\ref{components}
further yields that $w\cdot \eta =\widehat{ \sum_{i=1}^m \lambda_i \eta_i}$ where the coefficients $\lambda_i$'s are all non-negative.
That is, $w\cdot\eta$ is in the convex hull of affine limit roots. Then it follows from Proposition~\ref{hull} that
$w\cdot\eta \in E(W)$, and consequently, $\eta \in E(W)$.
\end{proof}

\begin{definition}
Define
\begin{align*}
E_{\afftype}&:=\{\eta\in E(W)\mid \#\pos(\eta) <\infty \text{\,and $\supp(w\cdot\eta)$ is }\\
& \text{disconnected for some $w\in W$ }\},\\
\noalign{\hbox{the set of limit roots that are non-trivial sums of affine limit roots.}}
\end{align*}		
\end{definition}

In a Coxeter system $(W, R)$ in which $W$ is an infinite irreducible non-affine Coxeter group and
$R$ is a finite generating set, it follows readily from Theorem~\ref{th:key}, Proposition~\ref{sum}, and
Proposition~\ref{hull} that
$E_{\nonaff}^2 \cap (E_{\aff}\uplus E_{\afftype}) =\emptyset$.
Furthermore, by combining Proposition~\ref{pp:im}, Theorem~\ref{th:key}, Proposition~\ref{sum}, and Theorem~\ref{ZQE} we immediately arrive at the following conclusion:
\begin{corollary}
\label{finite}
$$\mathscr{Z}\cap \widehat{Q}= E(W)\cap \mathscr{Z}=E_{\aff}\uplus E_{\afftype}=\{\eta\in E(W)\mid \#\pos(\eta)<\infty\}.$$
\qed
\end{corollary}

\section{Connection with dominance chains}

We close this paper with an application of the Pos-finite analysis developed in the previous section to the dominance partial ordering. In \cite[Open Question 7.7]{DHR13} it was asked: if $\{x_i\}_{i\in \N}\subset \Phi^+$ is an infinite injective sequence of positive roots totally ordered by dominance, that is, $x_{i+1}\SD x_i$ for all $i\in \N$, will the corresponding sequence of normalized roots $\{\widehat{x_i}\}_{i\in \N}$ possess a unique limit root. It has been observed in \cite[Propostion 6.8]{DHR13} that if such a normalized sequence of an infinite injective sequence of positive roots totally ordered by dominance possess two distinct accumulation points $\eta_1$ and $\eta_2$ then $(\eta_1, \eta_2)=0$ and they span a two dimensional totally isotropic subspace of $V$. In \cite[Lemma 2.10]{HMN} it has been further shown that in all Coxeter groups where the bilinear forms afforded by the Tits realization are of signature $(n-1, 1)$ such an infinite sequence possesses a unique accumulation point. Moreover in a forthcoming paper we show that in all infinite Coxeter groups of ranks $3$ and $4$, there is also unique accumulation. Based on these findings we have some confidence that the uniqueness of accumulation can be generally expected for all infinite Coxeter groups of finite rank. Although at this present stage we are unable to resolve this quesion in its entirety, in this section nevertheless we prove that if such a normalized sequence possesses an accumulation point $\eta\in E_{\aff}$ then this normalized sequence is convergent with $\eta$ being the only accumulation.

We believe that the significance of resolving this possible uniqueness of accumulation among infinite dominance chains is not to be underestimated.   So far, in the exploration of limit roots, examples of convergent infinite sequence of normalized roots are scarce, and indeed, for the time being we can only rely on infinite sequence coming from infinite dihedral reflection subgroups or affine reflection subgroups. Needless to say that such examples are severely limited in their scopes, and the need to construct other possibly more non-trivial examples of convergent infinite sequence of normalized roots is somewhat pressing. If one could prove that in a finitely generated infinite Coxeter group, the normalized sequence corresponding to an infinite injective sequence of normalized roots is convergent then we would have succeeded in someway of constructing non-trivial example of convergence. Furthermore, since dominance partial ordering is compatible with reflection subgroups (containing the reflections corresponding to those roots satisfying a dominance condition) such convergence behaviours work well with reflection subgroups.

We first list a number of observations concerning dominance, establishing certain inequalities which may prove to be useful, also we prove a result establishing when three distinct roots satisfying dominance  whose corresponding reflections generate an affine reflection subgroup. We recall the following technical result taken from \cite{FU2} (which is a refinement in one direction of \cite[Proposition 4.10]{FU2}  that for two roots $x$ and $y$, we have $x\dom y$ if and only if $x-y\in \mathscr{Z}\subset U^*$).
\begin{proposition}
 \label{pp:key}
Suppose that $x, y\in \Phi$ such that $x\SD y$. Then there exists
some $w\in W$ such that $wx\in \Phi^+$, $wy\in \Phi^-$ and $(w(x-y), z)\leq 0$
for all $z\in \Phi^+$.
\qed
\end{proposition}

From the above we may deduce a number of inequalities.
\begin{lemma}
\label{lem3}
Suppose that $a, b, c\in \Phi$ such that $a\dom b \dom c$.
Then $(a, c)\geq (b, c)$.
\end{lemma}
\begin{proof}
There is nothing to prove if either $a=b$ or $b=c$. Thus we may assume that the dominance
relations are strict. Since $a$ strictly dominates $b$, Proposition~\ref{pp:key} then yields that there
exists some $w\in W$ with $w a\in \Phi^+$, $w b\in \Phi^-$, and $(w(a-b), z)\leq 0$ for all $z\in \Phi^+$.
Since $b \dom c$ and $wb\in \Phi^-$, it follows that $wc\in \Phi^-$,  and consequently  $(w(a-b), wc)\geq 0$; that is
$(a-b, c)\geq 0$.
\end{proof}

An entirely similar proof as the above also establishes the following:
\begin{lemma}
\label{lem3'}
Suppose that $a, b, c\in \Phi$ such that $a\dom b \dom c$.
Then $(a, b)\leq (a, c)$.
\qed
\end{lemma}

\begin{proposition}
\label{eq1}
Suppose that $x, y, z$ are roots with $x\SD y\SD z$. If further $(x, z)=(x, y)$ then
$$(x, z)=(x, y)=(y, z)=1,$$
and the reflection subgroup $W':=\langle \rho_x, \rho_y, \rho_z\rangle$ is affine dihedral.
\end{proposition}
\begin{proof}
First, since $x\SD y\SD z$, it follows from  \cite[Proposition 4.10]{FU2} that $x-y\in  U^*$ and $y-z\in U^*$.
Then
$$0\geq (x-y, y-z)=\underbrace{(x, y-z)}_{=(x, y)-(x,z)=0}-\overbrace{(y, y-z)}^{=1-(y, z)\leq 0},$$
where $0\geq (x-y, y-z)$ follows from \cite[Theorem 5.8]{FU3}.
Consequently, $(y, y-z)=0$, and so $(y, z)=1$. Then
$$(x, y-z)=0, \quad (y, y-z)=0, \quad\text{and}\quad(z, y-z)=0.$$
That is $y-z$ is in the radical of the bilinear form $(\,,\,)$ restricted to the subspace spanned by $x$, $y$, and $z$.
Furthermore $$w\cdot (\widehat{y-z})=\widehat{y-z}, \quad \forall w\in W',$$
which in turn establishes that $E(W')=\{\widehat{y-z}\}$ is a singleton (following \cite[Theoerem 3.1]{DHR13}).
Now suppose for a contradiction that $W'$ is non-affine. Then it follows from Proposition~\ref{prop:aff} that there exists s non-affine dihedral reflection subgroup $W''\subset W'$ with $|E(W'')|=2$, contradicting $E(W'')\subseteq E(W')$. Consequently, $W'$ is affine, whence
$$(x, y)=(y, z)=(x, z)=1.$$
Then $(x-y, x)=(x-y, y)=(x-y, z)=0$, and $x-y$ is also in the radical of the bilinear form $(\,,\,)$ restricted to the subspace spanned by $x$, $y$, and $z$. Now if $W'$ is of rank $3$, then $\{x, y, z\}$ is linearly independent, and the radical of the bilinear form $(\,,\,)$ restricted to the subspace spanned by $x$, $y$, and $z$ is of dimension $2$ spanned by $x-y$ and $y-z$. On the other hand, it is well known (for instance, see Section 2.6 of \cite{HM}) that the radical of the restriction of the bilinear form $(\,,\,)$ on the subspace spanned by $x$, $y$, and $z$ is $1$-dimensional if $W'$ is irreducible affine, and we have a contradiction, and whence $W'$ is dihedral.

\end{proof}

\begin{proposition}
\label{domlimit}
Let $\{x_i\}_{i\in \N}\subset\Phi^+$ be an infinite sequence of positive roots satisfying $x_{i+1}\SD x_i$  for all $i\in \N$. Suppose that there is an  accumulation point $\eta$ of the corresponding normalized sequence $\{\widehat{x_i}\}_{i\in N}$ and $\eta\in E_{\aff}$. Then the normalized sequence $\{\widehat{x_i}\}_{i\in \N}$  is convergent with $\eta$ being the limit.
\end{proposition}
\begin{proof}
Suppose for a contradiction that there is yet another accumulation point for the normalized sequence $\{\widehat{x_i}\}_{i\in \N}$, say $\eta'$.
Choose two subsequences $\{y_i\}_{i\in \N}$ and $\{z_i\}_{i\in \N}$ of $\{x_i\}_{i\in \N}$ such that $\lim_{i\to\infty}\widehat{y_i}=\eta$ and $\lim_{i\to \infty}\widehat{z_i}=\eta'$, furthermore, we require that
$$z_{i+1}\SD y_{i+1}\SD z_i, \qquad\forall i\in \N.$$
It follows from Corollary~\ref{co:std} that there exists some $w\in W$ such that $\{w\cdot \eta\}=E(W_I)$ for some irreducible affine standard parabolic subgroup $W_I$. Note that $\{wx_i\}_{i\geq N}$ for some fixed $N\geq \ell(w)$ is still an infinite injective sequence of positive roots totally ordered by dominance; further, $w\cdot \eta$ and $w\cdot \eta'$ are now accumulation points for $\{\widehat{wy_i}\}_{i\geq N}$ and $\{\widehat{wz_i}\}_{i\geq N}$, respectively.
Since now $\pos(w\cdot \eta)=\emptyset$, it follows that
$$(w\cdot\eta, \widehat{ wy_i})\leq 0 \quad\text{and}\quad (w\cdot \eta, \widehat{ w z_i})\leq 0, \quad\forall i\geq N.$$
On the other hand, the dominance in the two subsequences implies that, $\forall i\geq N$,
\begin{align*}
(w\cdot \eta, \widehat{w y_i})&=\lim_{j\to \infty}(\widehat{ w y_j}\, ,\, \widehat{w y_i})\geq 0,\\
\noalign{\hbox{and}}
(w\cdot \eta, \widehat{w z_i})&=\lim_{j\to \infty}(\widehat{ w y_j}\, ,\, \widehat{w z_i})\geq 0.
\end{align*}
Thus, $\forall i\geq N$
$$(w\cdot \eta, \widehat{w y_i})=(w\cdot \eta, \widehat{w z_i})=0.$$

Note now that $I$, being a connected component of $\Pi$, is the set of simple roots for $W_I$, and furthermore, $I=\supp(w\cdot \eta)$.
The latter assertion is well-known and it can be seen in the following way: if $\supp(w\cdot\eta) \neq I$ then there exists some $a\in I$ while $a\notin\supp(w\cdot \eta)$, and the connectedness of $I$ implies that there must be some $b\in \supp(w\cdot \eta)$ such that $(a, b)<0$, but then $(w\cdot \eta, b)\neq 0$, contradicting that $w\cdot \eta\in \rad(I)$.
Now suppose that $x\in \Phi^+$ such that $(x, w\cdot \eta)=0$. Then the fact that $\supp(w\cdot\eta)=I$ readily implies that there is no  $b\in\supp(x)$ for which $(b, a)\neq 0$ for some $a\in I$. Thus, all $\supp(w z_i)$, $\forall i\geq N$, must either be disconnected from $I$ or contained in $I$ (since the support of a root must be itself connected). By our construction, $\supp(w y_i)= I$ for all $i$ large enough. Now if some $wz_j \notin \Phi(W_I)$ where $j$ is large enough, then $\supp(wz_j)$ must be disconnected from $I$, but this contradicts that $w y_j\,\SD\, w z_j$. Thus
$$wz_i\in \Phi(W_I), \quad \text{for all $i$ large enough}.$$
Since $E(W_I)=\{w\cdot \eta\}$, it follows that
$$w\cdot \eta'=\lim_{i\to \infty} \widehat{w z_i} =w\cdot \eta,$$
whence $\eta'=\eta$, and such a sequence has a unique limit.
\end{proof}

In a separate forthcoming paper, we obtain additional results establishing that the support of any accumulation point of a normalized sequence corresponding to an infinite injective sequence of roots totally ordered by dominance is connected (and remains connected upon the dot-action by the group, and therefore can not be of affine-type); further, we show that any two potential accumulation points $\eta_1$ and $\eta_2$ from the same normalized sequence corresponding to one infinite injective dominance chain have the same support, and $\pos(\eta_1)=\pos(\eta_2)$.

\section{Acknowledgments}

The first author is supported by the NSFC Grant (K19281004), and the third author is supported by the Research Development Fund from XJTLU (RDF-22-01-065), and they wish to sincerely acknowledge these support.

\bibliographystyle{amsplain}

\providecommand{\bysame}{\leavevmode\hbox to3em{\hrulefill}\thinspace}
\providecommand{\MR}{\relax\ifhmode\unskip\space\fi MR }
\providecommand{\MRhref}[2]{%
  \href{http://www.ams.org/mathscinet-getitem?mr=#1}{#2}
}
\providecommand{\href}[2]{#2}

\end{document}